\documentclass[reqno]{amsart}

\usepackage{amssymb,stmaryrd,amsthm,amsmath,bussproofs}
\usepackage[all]{xy}
\usepackage[latin1]{inputenc}
\usepackage{graphics}
\usepackage{float}

\newcommand{\myemph}[1]{\textbf{#1}}    

\newcommand{\type}{\textnormal{\texttt{type}}}       

\newcommand{\app}{\textnormal{\texttt{app}}} 
\newcommand{\pair}{\textnormal{\texttt{pair}}} 
\newcommand{\Jelim}{\textnormal{\texttt{J}}}
\newcommand{\id}[1]{\textnormal{\texttt{Id}}_{#1}} 
\newcommand{\idn}[2]{\underline{#1}^{#2}}
\newcommand{\judge}[3][]{#2\;\vdash_{#1}\;#3}

\newcommand{\globfor}[1]{\mathbf{\Gamma}(#1)}
\newcommand{\namefor}[1]{\ulcorner #1\urcorner}
\newcommand{\NN}{\mathtt{N}}
\newcommand{\zero}{\mathtt{0}}
\newcommand{\rr}{\mathtt{r}}
\newcommand{\successor}{\mathtt{S}}
\newcommand{\rec}{\textnormal{\texttt{rec}}}
\newcommand{\RSigma}{\textnormal{\texttt{R}}}

\newcommand{\DD}{\mathcal{D}}

\newcommand{\TT}{\mathcal{T}}
\newcommand{\sets}{\textnormal{\textbf{Set}}}

\newcommand{\groupoids}{\textnormal{\textbf{Gpd}}}

\newcommand{\gsets}{\textnormal{\textbf{GSet}}}
\newcommand{\rglob}{\textnormal{\textbf{rGSet}}}
\newcommand{\rgrph}{\textnormal{\textbf{rGraph}}}
\newcommand{\rgsets}{\textnormal{\textbf{rGSet}}}
\newcommand{\mlcxs}{\textnormal{\textbf{MLCx}}}
\newcommand{\extomega}{\textnormal{\textbf{Ext}}(\mathbb{T}_{\omega})}
\newcommand{\monadsin}{\textnormal{\textbf{Mon}}}
\newcommand{\freegpd}{\mathcal{F}}
\newcommand{\CC}{\mathcal{C}}
\newcommand{\II}{\mathcal{I}}

\newtheorem{theorem}{Theorem}[section]
\newtheorem{lemma}[theorem]{Lemma}
\newtheorem{proposition}[theorem]{Proposition}
\newtheorem{corollary}[theorem]{Corollary}

\theoremstyle{definition}
\newtheorem{definition}[theorem]{Definition}
\newtheorem{example}[theorem]{Example}

\theoremstyle{remark}
\newtheorem*{remark}{Remark}

\newcommand{\colimit}{\varinjlim}

\newcommand{\myexp}{\operatorname{Exp}}

\newcommand{\lscott}{[\![}
\newcommand{\rscott}{]\!]}

\newcommand{\iso}{\cong}
\newcommand{\homotopic}{\simeq}
\newcommand{\adjoint}{\dashv}
\newcommand{\forces}{\Vdash}

\floatstyle{ruled}
\newfloat{mytable}{htbp}{lot}[section]
\floatname{mytable}{Figure}

\begin{document}

\title{martin-l\"{o}f complexes}
\author{S. Awodey}
\address{Department of Philosophy, Carnegie Mellon University\\5000
  Forbes Ave., Pittsburgh, PA 15213 USA}
\email{awodey@cmu.edu}
\thanks{Awodey gratefully acknowledges the support of the National
  Science Foundation, Grant DMS-1001191 and the Air Force OSR, Grant
  11NL035.}

\author{P. Hofstra}
\address{Department of Mathematics and Statistics, University of
  Ottawa\\585 King Edward Ave., Ottawa, ON K1N 6N5 Canada}
\thanks{Hofstra is partially supported by NSERC, Canada.}
\email{phofstra@uottawa.ca}

\author{M.~A. Warren}
\address{School of Mathematics, Institute for Advanced Study\\Einstein Drive, Princeton, NJ 08540 USA}
\thanks{During the preparation of this work Warren received support
  from the Fields Institute and the Atlantic Association for Research in
  the Mathematical Sciences.  He is currently supported by the National Science
  Foundation.  In particular, this material is based upon work
  supported by the National Science Foundation under agreement
  No. DMS-0635607. Any opinions, findings and conclusions or
  recommendations expressed in this material are those of the authors
  and do not necessarily reflect the views of the National Science Foundation.}
\email{mwarren@math.ias.edu}

\dedicatory{Dedicated to Per Martin-L\"{o}f on the occasion of his retirement.}

\date{\today}

\begin{abstract}
  In this paper we define Martin-L\"{o}f complexes to be
  algebras for monads on the category of (reflexive) globular sets
  which freely add cells in accordance with the rules of intensional
  Martin-L\"{o}f type theory.  We then study the resulting categories
  of algebras for several theories.  Our principal result is that there exists a cofibrantly generated Quillen model
  structure on the category of 1-truncated Martin-L\"{o}f complexes
  and that this category is Quillen equivalent to the category of
  groupoids.  In particular, 1-truncated Martin-L\"{o}f complexes are
  a model of homotopy 1-types.
\end{abstract}

\maketitle

\tableofcontents

\section{Introduction}\label{sec:intro}
This paper pursues a surprising connection between Geometry, Algebra,
and Logic that has only recently come to light, in the form of an
interpretation of the constructive type theory of Martin-L\"of into
homotopy theory, resulting in new examples of certain algebraic
structures which are important in topology.  This fascinating
connection is currently under investigation from several different
perspectives (\cite{Awodey:HMIT,Warren:PhD,Garner:TDMTT,Lumsdaine:WOITT,Berg:TWOG,Gambino:ITWFS}
), and these preliminary results confirm the significance of the link.  Some of these results will be surveyed in this brief introduction in order to position the present work in its context; especially for the reader coming from one field or the other, a brief summary is given of the essential concepts from the different subjects involved.

Martin-L\"of type theory \cite{MartinLof:ITT} is a formal system originally intended to provide a rigorous
framework in which to develop constructive mathematics. At heart, it is a calculus
for reasoning about dependent types and terms, and equality between those.  Under the
Curry-Howard correspondence, one may identify types with propositions, and terms with proofs.
Viewed in this manner, the system can be shown to be at least as strong as second-order logic,
and it is also known to interpret constructive set theory. Indeed, Martin-L\"of type theory has 
been used successfully to formalize parts of constructive mathematics, such as
pointless topology (constructive locale theory). Moreover, it has been employed as a framework
for the development of programming languages as well, a task for which it is especially well-suited in virtue of its combination of expressive strength and desirable proof-theoretic properties. (See the textbook~\cite{Nordstrom:PMLTT} for a discussion.)

The type theory has two variants: an intensional, and an extensional version. The difference between them lies mainly in the treatment of equality. In the intensional version (with which we are mainly concerned in the present work),  one has two different kinds of equality: the first kind is called \myemph{definitional equality}, and behaves much like equality between terms in the simply-typed lambda-calculus, or any other conventional equational theory. The second kind is a more subtle relation, called \myemph{propositional equality}, which, under the Curry-Howard  correspondence, represents the equality formulas of first-order logic. Specifically, given two terms $a,b$\/ of the same type $A$, one may form a new type $\id{A}(a,b)$, which we think of as the proposition that $a$\/ and $b$\/ 
are equal; a term of this type thus represents a proof of the proposition that $a$ equals $b$ (hence the name ``propositional equality''). 

When $a$\/ and $b$\/ are definitionally equal, then (since they can be freely substituted for each other) they are also propositionally equal, in the sense that the type $\id{A}(a,b)$ is inhabited by a term; but the converse is generally not true, at least in the intensional version of the theory.  In the extensional version, by contrast, the two notions of equality are forced by an additional rule to coincide. As is well-known, however, the price one pays for this simplification is a loss of desirable proof-theoretic properties, such as strong normalization
and decidable equality of terms.

In the intensional version with which we shall be concerned here, it
can be shown that the identity types $\id{A}(a,b)$\/ carry certain
structure which was observed by Hofmann and Streicher in
\cite{Hofmann:GITT} to be analogous to that of a groupoid. Specifically, the reflexivity of propositional equality produces identity proofs $\rr(a):\id{A}(a,a)$\/ for any term $a:A$, playing the role of a unit arrow for $a$; and when $f:\id{A}(a,b)$\/ is an identity proof, then (corresponding to the symmetry of identity) there also exists a proof $f^{-1}:\id{A}(b,a)$, to be thought of as the inverse of $f$;
finally, when $f:\id{A}(a,b)$\/ and $g:\id{A}(b,c)$\/ are identity proofs, then (corresponding to transitivity) there is a new proof $(g \cdot f):\id{A}(a,c)$, thought of
as the composite of $f$\/ and $g$. Moreover, this structure on each type $A$ can be shown to satisfy the usual groupoid laws, but only \myemph{up to propositional equality}.  We shall return to this point below.

\subsection{Groupoid semantics}

A good notion of a model for the extensional theory is due to Seely~\cite{Seely:LCCCTT}, who showed that one can interpret type dependency in locally cartesian closed categories in a very natural way. (There are certain coherence issues related to this semantics, prompting a later refinement, but this need not concern us here.) Of course, intensional type theory may also be interpreted this way in lcccs, but then the interpretation of the identity types necessarily becomes trivial.

The first non-trivial semantics for intensional type theory was
developed by Hofmann and Streicher~\cite{Hofmann:GITT} using
\myemph{groupoids}, which are categories in which every arrow is an
isomorphism.   The category of groupoids is not locally cartesian
closed, and the model employs certain fibrations (equivalently,
groupoid-valued functors) to model type dependency.  A closed type $A$\/ will be 
interpreted as a groupoid, a term $a:A$\/ as an object of this groupoid, and an identity proof
$f:\id{A}(a,b)$\/ as an arrow $f:a\to b$ in $G$.  The
interpretation no longer validates extensionality, since there can be
different elements $a, b$ related by non-identity arrows $f:a\to b$.
Indeed, there may be many different such arrows $f,g,\dots: a\to b$\,;
however, unlike in the type theory, these cannot in turn be further
related by non-trivial identity terms of higher type $\vartheta:\id{\id{A}}(f,g)$,
since a (conventional) groupoid has no non-trivial higher-dimensional structure.  Thus the groupoid semantics validates a certain truncation principle, stating that all higher identity types are trivial---a form of extensionality one dimension up.
In particular, the groupoid laws for the identity types are strictly satisfied in these models, rather than holding only up to propositional equality.

This situation has led to the use of the higher-dimensional analogues
of groupoids, as formulated in category theory, in order to provide
models admitting non-trivial higher identity types.   Such higher
groupoids have been studied extensively in homotopy theory in recent
years, since they occur naturally as the (higher) fundamental
groupoids of spaces (see below).   In this direction, Warren
\cite{Warren:PhD} has generalized the groupoid model of
\cite{Hofmann:GITT} to strict $\omega$-groupoids, thereby showing that
the type theory truly possesses non-trivial higher-dimensional
structure.  Along similar lines, Garner~\cite{Garner:TDMTT} has used a
2-dimensional notion of fibration to model intensional type theory,
and shown that when various truncation axioms are added the theory is
sound and complete with respect to this semantics.

\subsection{Homotopy theory}

In homotopy theory one is concerned with spaces and continuous mappings up to homotopy; a homotopy between continuous maps $f,g:X \to Y$\/ is a continuous map $\vartheta:X \times [0,1] \to Y$\/ satisfying $\vartheta(x,0)=f(x)$\/ and $\vartheta(x,1)=g(x)$. Such a homotopy $\vartheta$ can be thought of as a ``continuous deformation" of $f$ into $g$, determining a higher-dimensional arrow $\vartheta: f \to g$.  As already suggested, one also considers homotopies between homotopies, referred to as higher homotopies. Algebraic invariants, such as homology or the fundamental group, are homotopy-invariant:  they are invariants of the \emph{homotopy types} of spaces, i.e. of
equivalence classes of spaces under the homotopy equivalence relation, where two spaces $X$\/ and $Y$\/ are said to be homotopy equivalent if there exist maps $f:X \to Y$\/ and $g:Y \to X$\/ for
which the composites $gf$\/ and $fg$\/ are homotopic to the identity maps on $X$\/ and on $Y$, respectively.

When we consider the points of a space $X$, the paths in $X$, the homotopies between
paths, and all higher homotopies, we  obtain a structure called the
\myemph{fundamental weak $\omega$-groupoid of $X$}. We can truncate
this structure by considering only the points and paths up to homotopy, and this
yields the usual fundamental groupoid of the space.  This truncation
is evidently analogous to adding to our type theory axioms of
extensionality above the first identity type. Indeed, these constructions, including the basic
assignment of fundamental groupoids to objects, are special cases of a
common, general construction that can be described abstractly in
axiomatic homotopy theory.  The central concept is that of a
\myemph{Quillen model category}, which captures axiomatically some of
the essential features of homotopy of topological spaces, enabling us
to ``do homotopy" in different mathematical settings, and to express
the fact that two categories carry the same homotopical information,
even if they are not equivalent in the ordinary sense.  The basic
result of Awodey and Warren in~\cite{Awodey:HMIT} (see
also~\cite{Warren:PhD}) is that it is possible to model the type
theory in any Quillen model category which is well-behaved in certain
ways (essentially using just the basic notion of a weak factorization
system).  In this interpretation, one uses path objects to model
identity types in a non-trivial way, recovering the groupoid model as
a special case. This suggests that intensional type theories are a
sort of internal language of (certain kinds of) model categories.
Indeed, in \cite{Gambino:ITWFS} it is shown that the type theory
itself carries a natural such homotopy structure (i.e.\ a weak
factorization system), so that the theory is not only sound but also
complete with respect to such abstract homotopical semantics.

Thus we are justified in thinking of types in the intensional theory as \myemph{spaces}. From this point of view, the
terms of the type $A$\/ are the points of the ``space" $A$, the identity type $\id{A}(a,b)$\/ represents the
collection of paths from $a$\/ to $b$, and the higher identities are homotopies between paths, 
homotopies between homotopies of paths, et cetera. The fact that paths and homotopies do not form a groupoid, but only a groupoid up to homotopy, is of course precisely the same observation as the fact that the identity types only satisfy the groupoid laws up to propositional equality. This parallel between type theory and
homotopy theory, which was first pointed out by Moerdijk a few years ago, has now been made precise by the recognition that both cases are instances of one and the same abstract axiomatic theory.

In particular, it has been shown independently by
Lumsdaine~\cite{Lumsdaine:WOITT} and Van den Berg and
Garner~\cite{Berg:TWOG} that the tower of identity types over any
fixed base type $A$\/ in the intensional theory indeed gives rise to a
certain infinite dimensional categorical structure called a weak
$\omega$-groupoid. In fact, something apparently stronger is shown,
namely that at every type the type theory already hosts an internal
model of such a higher category.  The next step in exploring the
connection between type theory and topology
is to investigate the relationship between type theoretic
``truncation" (i.e.\ higher-dimensional extensionality principles) and
topological ``truncation" of the higher fundamental groups. Spaces for
which the homotopy type is already completely determined by the
fundamental groupoid are called \myemph{homotopy 1-types}, or simply
1-types. More generally, one has $n$-types, which are thought of as
spaces which have no homotopical information above dimension $n$. One of
the goals of homotopy theory is to obtain good models of homotopy
$n$-types.  For example, the category of groupoids is Quillen equivalent
to the category of 1-types and therefore the corresponding homotopy
categories (obtained by inverting weak equivalences) are equivalent;
in this precise sense, groupoids are said to model homotopy 1-types
(for more on homotopy types see \cite{Baues:HT}).  

\subsection{Contributions of this paper}

The current paper aims at further investigation of the relationship
between type theory and homotopy theory, but in a way that is somewhat
different from the work already mentioned. First of all, our primary objective
is not to give a new semantics, although some of the results
will depend on a new model which will be presented in a sequel to this paper. Secondly,
while earlier work centered around constructing
higher-dimensional structures from type theories, we are also
interested in understanding the limitations of this process. Finally,
we wish to make another connection between model categories and type
theory, namely by showing that a category of suitably truncated type
theories gives a model of the homotopy 1-types. It is our hope that
this picture can then be extended to higher dimensions.

Our first goal is to show how every extension of intensional type theory gives rise to a monad on the category of globular sets.
Intuitively, the monad associated to a theory freely adds cells to a
globular set in accordance with the structure imposed on the tower of identity types over a base type by the rules of the type theory.
For example, the monad will formally add composites and inverses for
all cells of dimension strictly greater than 0 in the globular set; however, it
adds much more than just these formal composites; it also produces a plethora of new cells which we here call \myemph{doppelg\"angers}.
For every such monad we may consider its category of algebras: these
we refer to as \myemph{Martin-L\"of complexes} (or \myemph{ML-complexes}), and these are
the main objects of study of the paper.

The theories which we shall consider arise from basic intensional Martin-L\"of type theory 
having dependent sums and products as well as a natural numbers object.
(The latter plays no conceptual role in this paper 
but because of its importance in virtually every application of the theory to mathematics and computer science we thought
it important to show that our results are not affected by its presence.)
We shall then consider extensions of this basic theory obtained by adding truncation axioms, which effectively trivialize the higher identity types
above a fixed dimension. Using these theories we get a hierarchy of categories of Martin-L\"of complexes, and in this paper we shall investigate the
first two dimensions in detail.

The 0-dimensional case is relatively straightforward --- we shall prove here that the monad on globular sets is idempotent and is
in fact isomorphic to the connected components functor, so that its
category of algebras (the 0-dimensional ML-complexes) is equivalent to the
category of sets.

Matters become more interesting in dimension 1. Towards an analysis of
1-dimensional ML-complexes we first observe, using the
Hofmann-Streicher groupoid semantics, that every ML-complex has an underlying groupoid, 
and that there is a canonical comparison functor between the
underlying groupoid of a free ML-complex and the free groupoid on the
same globular set.  This functor is not an isomorphism of groupoids,
because the free ML-complex is, intuitively speaking, much larger due
to all the doppelg\"angers produced by the theory. The main technical
difficulty then is to prove that there is still an equivalence of
groupoids between the two. This result follows from a proof-theoretic analysis of the theories in question, 
which will be borrowed from~\cite{Hofstra:CRMTT}. One of the main results in loc. cit. 
allows us to conclude that every term of the theory represents, up to propositional equality, an object or morphism of the free groupoid. 
This essentially shows that even though the theory forces the existence of many more objects and arrows than needed to form the free
groupoid, it does not force anything which is undesirable from a homotopical point of view. 

Once this key result is in place, we turn to an analysis of the category of 1-dimensional ML-complexes as a whole. To start, we set up an adjunction
between this category and the category of groupoids. This adjunction is not an equivalence: a ML-complex structure on a globular set carries
essentially more information than a groupoid structure. We can, however, make use of the adjunction by transferring along it the standard Quillen model
structure on the category of groupoids \cite{Joyal:SSCS}, turning the category of 1-dimensional ML-complexes into a cofibrantly generated model category. 

Finally, we prove that the adjunction between groupoids and 1-dimensional ML-complexes is in fact a Quillen equivalence. Because the categories
of groupoids and that of homotopy 1-types are Quillen equivalent, this
makes precise in which sense the 1-truncated version of the type theory models homotopy 1-types. 
It also explains why the groupoid semantics is adequate from a homotopical point of view, but is still incomplete because it lacks the
possibility (which is present in ML-complexes) of handling different interpretations for doppelg\"anger terms.

\subsection{Plan of the paper}

In Section \ref{sec:background} we recall the basics of Martin-L\"{o}f
type theory as well as several facts about globular sets which will be
required later.  We also fix notation (some of which is
non-standard).  The reader who is familiar with this material should
feel free to skip ahead. 
 
Section \ref{sec:monad} describes the
construction of monads on the category of reflexive globular sets
coming from type theories.  We then define the categories $\mlcxs$ and
$\mlcxs_{n}$ of Martin-L\"{o}f complexes and $n$-truncated
Martin-L\"{o}f complexes as the Eilenberg-Moore categories of the
monads $M_{\omega}$ and $M_{n}$, respectively, generated by suitable
theories.  These monads are shown to be finitary and it therefore
follows that the categories $\mlcxs$ and $\mlcxs_{n}$ are complete and
cocomplete.

In Section
\ref{sec:doppelgangers} we study $0$-truncated and $1$-truncated Martin-L\"{o}f
complexes.  We first show that the category
$\mlcxs_{0}$ is equivalent to the category of sets and, moreover, that
if $X$ is a reflexive globular set, then $M_{0}(X)$ is the set of
connected components of $X$.  
Even the proofs of these eminently plausible results
are a bit more complicated than one might at first expect;  
one of the principal difficulties
one faces when proving results about Martin-L\"{o}f complexes is that
the type theory also adds, in addition to composition and inverses,
the doppelg\"{a}nger terms mentioned earlier. 
We then turn to $1$-truncated complexes, 
with a proof that every such complex can be
equipped with the structure of a groupoid. Towards a characterization of
the free 1-dimensional ML-complexes, it is first shown that
the Hofmann-Streicher groupoid semantics induces a comparison functor
between the free groupoid $\freegpd(G)$ on a reflexive globular set 
and the induced groupoid structure on $M_{1}(G)$, the free 1-dimensional algebra on $G$. 
The main technical observation, namely that this comparison functor is 
an equivalence of groupoids, follows from the realizability semantics
presented in~\cite{Hofstra:CRMTT}. Because this technique will also be required
later, we shall give, for reasons of self-containment of
the present paper, a brief explanation of how it works.

Finally, Section~\ref{sec:model} shows that $\mlcxs_{1}$\/ can be endowed 
with a Quillen closed model structure.  This model structure is obtained from
an adjunction with the category of groupoids via Quillen's path object
argument \cite{Quillen:HA}. The main result states that the adjunction
between $\mlcxs_{1}$ and the category of groupoids is a Quillen
equivalence (Theorem \ref{thm:equiv}).  Both the
existence of the model structure and the verification of the Quillen
equivalence make crucial use of the realizability semantics from \cite{Hofstra:CRMTT}.

\subsection*{Acknowledgements}

We are greatly indebted to the anonymous referee for providing us with
numerous insightful comments and suggestions and for pointing out, as
Peter LeFanu Lumsdaine also did, a gap in an earlier version of this
paper. We would also like to thank Nicola Gambino,
Richard Garner, Peter LeFanu Lumsdaine, Thomas Streicher and Phil Scott for
useful discussions of some of the ideas in this paper.
\section{Background}\label{sec:background}

The purpose of this section is to provide the reader with a brief
introduction to Martin-L\"of type theory. We begin by giving a quick exposition of
the main features of the most basic version of the theory we shall be
concerned with. In particular we explain the different kinds of judgements of
the system, dependent products and sums, identity types and the notion of propositional equality. We also use this as an opportunity
to fix some notation and terminology, in particular concerning identity types.  

We assume that the reader is somewhat familiar with at least simple type theory. For more background on (dependent) type theory we refer to the
textbook~\cite{Jacobs:CLTT}.  The reader who is more familiar with
higher-dimensional category theory or homotopy theory might also
consult \cite{Awodey:HMIT} for a ``homotopical'' view of type theory.

In the last subsection we introduce the basic categorical structures used in the paper, namely globular sets. A more detailed exposition of
globular sets may be found in~\cite{Street:PTGS}, or the textbook~\cite{Leinster:HOHC}.

\subsection{Type dependency, contexts and judgements}\label{subsec:judgements}

Type dependency means that types may depend on variables of other types; for example one can has a type $T(x)$\/ depending on a variable $x$\/ of type $S$.
Such a type $T(x)$\/ is often thought of as being indexed by the type $S$. To illustrate this, suppose that we let $S$\/ denote the type of rings; then the
type $T(x)$\/ of modules depends on, or varies over, the type of rings.  

One may then substitute a term $a$\/ of type $S$\/ into the type $T(x)$, as to obtain a new type $T(a)$. In the above example, $T(a)$\/ would be the
type of modules over the ring $a$.

The fact that types may depend on terms has two obvious consequences: first, one can no longer, as in simple type theory, separate the formation of
types and that of terms into two inductive defitions; rather, types and terms are derived simultaneously. Second, the notion of a \myemph{variable
context} also needs to take dependency into account. Explicitly, this means that a variable context $\Gamma$\/ is now an ordered sequence of variable
declarations $\Gamma=(x_1:T_1, \ldots, x_k:T_k)$, where each type $T_i$\/ may only depend on the variables declared earlier, i.e. on $x_1, \ldots, x_{i-1}$.
For example,
\[ x:S, y:S, z:T(x), v:R(x,y,z) \]
is a legitimate variable context, but
\[ x:S, y:s, v:R(x,y,z), z:T(x) \]
is not, because of the fact that $R$\/ depends on $z$, which hasn't been declared yet.
Throughout, we shall always assume that contexts are well-formed in this sense.

Thus the theory is concerned with types and terms in context, and with equalities between such types and terms. Formally,
statements about these are called \myemph{judgements}, and these come in four kinds:
\[ \judge{\Gamma}{T:\type}\]
This judgement states that $T$\/ is a type, possibly depending on the variables declared in the context $\Gamma$.
\[ \judge{\Gamma}{\tau:T}\]
This judgement states that $\tau$\/ is a term of type $T$, where both $\tau$\/ and the type $T$\/ may depend on the variables from $\Gamma$.
\[ \judge{\Gamma}{T=S:\type}\]
This judgement states that $T$\/ and $S$\/ are (definitionally) equal types.
\[ \judge{\Gamma}{\tau=\tau':T}\]
This judgement states that $\tau$\/ and $\tau'$\/ are (definitionally) equal terms of type $T$.

In the theory, such judgements are derived from axioms using inference rules. These derivations (which may formally be regarded as finite trees suitably labelled by
judgements and inference rules) are the main objects of study. Below we shall discuss several of the rules which may be used to derive new judgements from old;
the axioms typically include judgements stating the existence of certain basic types and terms. 

When the context plays no role in a judgement or rule of the theory, we shall usually omit it altogether.

\subsection{Definitional equality}\label{subsec:definitional}
The notion of equality here is the standard one, but often the qualifier \myemph{definitional} is used to distinguish it from the different
notion of propositional equality, to be discussed below. The rules governing the behaviour of definitional equality are as expected. 
Apart from the rules expressing that definitional equality is an equivalence relation, there are rules which force that it is a congruence with respect to
substitution into types and terms:
\begin{prooftree}
  \AxiomC{$\judge{}{a=b:A}$}
  \AxiomC{$\judge{x:A}{B(x):\type}$}
  \BinaryInfC{$\judge{}{B(a)=B(b):\type}$}
  \end{prooftree}
\begin{prooftree}
  \AxiomC{$\judge{}{a=b:A}$}
  \AxiomC{$\judge{x:A}{f(x):B(x)}$}
  \BinaryInfC{$\judge{}{f(a)=f(b):B(a)}$}
\end{prooftree}
 \begin{prooftree}
  \AxiomC{$\judge{}{A=B:\type}$}
  \AxiomC{$\judge{}{a:A}$}
  \BinaryInfC{$\judge{}{a:B}$}
\end{prooftree}
The first rule states that substituting equal terms into a type results in equal types; the second states the same, but now for substitution into terms;
the last rule states that equal types are inhabited by the same terms.
A complete set of rules for definitional equality may be found in the appendix.

\subsection{Dependent products and sums}\label{subsec:dependent}

There are several ways to construct new types from old. For each new type one specifies three things: an introduction rule which generates new terms
of the type; an elimination rule which shows how general terms of the new type may be used; and a conversion rule which governs the interaction between
the two.  

We now discuss the formation of dependent products and sums. Given a type $B(x)$\/ depending on 
a variable $x$\/ of type $A$, we may form the type
$\prod_{x:A}B(x)$, to be thought of as the type of sections of $B(x)$\/ over $A$. The rules are as follows:

\begin{prooftree}
  \AxiomC{$\judge{x:A}{B(x):\type}$}
  \RightLabel{$\prod$ formation}
  \UnaryInfC{$\judge{}{\prod_{x:A}B(x):\type}$}
\end{prooftree}
\vspace{.25cm}
\begin{prooftree}
  \AxiomC{$\judge{x:A}{f(x):B(x)}$}
  \RightLabel{$\prod$ introduction}
  \UnaryInfC{$\judge{}{\lambda_{x:A}.f(x):\prod_{x:A}B(x)}$}
\end{prooftree}
\vspace{.5cm}
\begin{prooftree}
  \AxiomC{$\judge{}{f:\prod_{x:A}B(x)}$}
  \AxiomC{$\judge{}{a:A}$}
  \RightLabel{$\prod$ elimination}
  \BinaryInfC{$\judge{}{\app(f,a):B(a)}$}
\end{prooftree}
\vspace{.25cm}
\begin{prooftree}
  \AxiomC{$\judge{}{\lambda_{x:A}.f(x):\prod_{x:A}B(x)}$}
  \AxiomC{$\judge{}{a:A}$}
  \RightLabel{$\prod$ conversion}
  \BinaryInfC{$\judge{}{\app\bigl(\lambda_{x:A}.f(x),a\bigr) \;=\; f(a):B(a)}$}
\end{prooftree}
Thus an introduction term of type $\prod_{x:A}B(x)$\/ is a lambda expression, thought of as an operation assigning to each $x:A$\/ a value $f(x):B(x)$.
A general term $f$\/ of type $\prod_{x:A}B(x)$\/ may be applied to a term $a$\/ of type $A$, as to return a term $\app(f,a)$\/ of type $B(a)$. Finally,
the conversion rule, commonly known as \myemph{beta-conversion}, allows us to reduce $\app(\lambda_{x:A}.f(x),a)$\/ to $f(a)$.
In the case where the type $B(x)$\/ does not depend on the variable $x$, we shall often write $B^A$\/ for the type $\prod_{x:A}B(x)$.

Similarly, the theory admits formation of \myemph{dependent sum types} $\sum_{x:A}B(x)$. The rules are:

\begin{prooftree}
  \AxiomC{$\judge{x:A}{B(x):\type}$}
  \RightLabel{$\sum$ formation}
  \UnaryInfC{$\judge{}{\sum_{x:A}B(x):\type}$}
\end{prooftree}
\vspace{.25cm}
\begin{prooftree}
  \AxiomC{$\judge{}{a:A}$}
  \AxiomC{$\judge{}{b:B(a)}$}
  \RightLabel{$\sum$ introduction}
  \BinaryInfC{$\judge{}{\pair(a,b):\sum_{x:A}B(x)}$}
\end{prooftree}
\vspace{.5cm}
\begin{prooftree}
  \AxiomC{$\judge{}{p:\sum_{x:A}B(x)}$}
  \AxiomC{$\judge{x:A,y:B(x)}{\psi(x,y):C\bigl(\pair(x,y)\bigr)}$}
  \RightLabel{$\sum$ elimination}
  \BinaryInfC{$\judge{}{\RSigma\bigl([x:A,y:B(x)]\psi(x,y),p\bigr):C(p)}$}
\end{prooftree}
\begin{prooftree}
  \AxiomC{$\judge{}{a:A}$}
  \AxiomC{$\judge{}{b:B(a)}$}
  \AxiomC{$\judge{x:A,y:B(x)}{\psi(x,y):C\bigl(\pair(x,y)\bigr)}$}
  \RightLabel{$\sum$ conversion}
  \TrinaryInfC{$\judge{}{\RSigma\bigl([x:A,y:B(x)]\psi(x,y),\pair(a,b)\bigr) \;=\; \psi(a,b):C\bigl(\pair(a,b)\bigr)}$}
\end{prooftree}
The notation $[x:A,y:B(x)]$ indicates that the variables $x$ and $y$
are formally bound in the term.  Using these rules, we may define projection terms by letting 
\[ \pi_0(p)=\RSigma([x:A,y:B(x)]x,p), \qquad \pi_1(p)=\RSigma([x:A,y:B(x)]y:B(p)\]  in

\begin{prooftree}
  \AxiomC{$\judge{}{p:\sum_{x:A}B(x)}$}
  \AxiomC{$\judge{x:A,y:B(x)}{x:A}$}
   \RightLabel{$\sum$ elimination}
  \BinaryInfC{$\judge{}{\pi_0(p):A}$}
\end{prooftree}

\begin{prooftree}
  \AxiomC{$\judge{}{p:\sum_{x:A}B(x)}$}
  \AxiomC{$\judge{x:A,y:B(x)}{y:B(x)}$}
   \RightLabel{$\sum$ elimination}
  \BinaryInfC{$\judge{}{\pi_1(p):B(\pi_0(p))}$}
\end{prooftree}
The projection terms $\pi_0(p)$\/ and $\pi_1(p)$\/ then
satisfy the conversion rules 
\begin{align*}
  \pi_{i}(\pair(a,b)) & = \begin{cases} 
    a & \text{ if }i=0\\
    b & \text{ if }i=1.
  \end{cases}
\end{align*}

We point out that we do not adopt the $\eta$-rule for sums 
\begin{prooftree}
 \AxiomC{$\judge{}{p:\sum_{x:A}B(x)}$}
 \UnaryInfC{$\judge{}{p=\pair(\pi_0(p),\pi_1(p)):\sum_{x:A}B(x)}$}
\end{prooftree}
but that it can easily be proved that every term of type 
$\sum_{x:A}B(x)$\/ is \emph{propositionally} equal to a pair term (see the discussion of identity types below for what this means). 

In some treatments a different formulation of the rules for sum types is used, taking the projection terms as primitive. In the presence of the
$\eta$-rule both formulations are equivalent, but without the
$\eta$-rule this latter approach is strictly weaker (see~\cite{Hofmann:SSDT}).

\subsection{Identity types}\label{subsec:identity}

Let $A$\/ be a type. For every pair of terms $a,b$\/ of type $A$\/ we may form
a new type $\idn{A}{}(a,b)$. This type is thought of as the type of proofs of the fact that $a$\/ and $b$\/ are equal. A term $\tau:\idn{A}{}(a,b)$\/ is
sometimes referred to as a \myemph{propositional identity proof}. It is important to note that the existence of such a proof term does not necessarily imply that
$a=b$\/ in the definitional sense of equality discussed above. 
From a more geometric perspective one may think of a propositional
equality as a \myemph{homotopy} between $a$\/ and $b$ (see~\cite{Awodey:HMIT}). This explains why we sometimes
use the notation $a \simeq b$\/ to indicate the existence of a propositional identity between $a$\/ and $b$. We also point out that it is perhaps more common to
denote the identity type $\idn{A}{}(a,b)$\/ by $\id{A}(a,b)$, but we have chosen to adopt a notation more suggestive of hom-sets.

The formation rule for the identity types is thus as follows (omitting contexts for simplicity)

\begin{prooftree}
  \AxiomC{$\judge{}{a,b:A}$}
  \RightLabel{$\id{}$ formation}
  \UnaryInfC{$\judge{}{\idn{A}{}(a,b):\type}$}
\end{prooftree}

\noindent where we write $a,b:A$ as an abreviation for the two judgements $a:A$
and $b:A$.  Then, there are the introduction and elimination rules:

\begin{prooftree}
  \AxiomC{$\judge{}{a:A}$}
  \RightLabel{$\id{}$ introduction}
  \UnaryInfC{$\rr(a):\idn{A}{}(a,a)$}
\end{prooftree}
\vspace{.2cm}
\begin{prooftree}
  \AxiomC{$\judge{x:A,y:A,z:\idn{A}{}(x,y)}{B(x,y,z):\type}$}
  \noLine
  \UnaryInfC{$\judge{x:A}{\varphi(x):B\bigl(x,x,\rr(x)\bigr)}$}
  \noLine
  \UnaryInfC{$\judge{}{f:\idn{A}{}(a,b)}$}
  \RightLabel{$\id{}$ elimination}
  \UnaryInfC{$\judge{}{\Jelim_{[x,y:A,z:\idn{A}{}(x,y)]B(x,y,z)}\bigl([x:A]\varphi(x),a,b,f):B(a,b,f)}$}
\end{prooftree}

\vspace{.25cm}
\noindent The introduction term $\rr(a)$\/ is called the \myemph{reflexivity term}; it witnesses the fact that $a \simeq a$. The elimination rule is a bit more involved.
What we start with is first of all a type $A$\/ (which is referred to as the \myemph{type over which the elimination occurs}), and an identity proof 
$f: \idn{A}{}(a,b)$\/ (this is the term which is being eliminated). Next we need a type $B(x,y,z)$\/ (called the \myemph{pattern type}) and a term $\varphi(x)$\/ 
of type $B(x,x,\rr(x))$\/ which intuitively witnesses the fact that the pattern type is inhabited in the trivial instance where we substitute a reflexivity term.
Given all of this, we may form a new term $\Jelim([x:A]\varphi(x),a,b,f)$\/ of type $B(a,b,f)$. One way to think of this $J$-term is as the result of
expanding the term $\phi(x)$\/ using the propositional equality $f:\idn{A}{}(a,b)$. This viewpoint will be developed in more detail later on.

Note also that the variables $x,y,z$\/ in the elimination rule need not necessarily occur in the type $B(x,y,z)$, and similarly that $x$\/ need not occur in
$\varphi(x)$. Also, it may happen that the term $f$\/ (and possibly also $a,b$) are themselves variable, in which case the $J$-term depends on those variables.

Finally, there is a conversion rule:

\begin{prooftree}
  \AxiomC{$\judge{}{a:A}$}
  \RightLabel{$\id{}$ conversion}
  \UnaryInfC{$\judge{}{\Jelim_{[x,y:A,z:\idn{A}{}(x,y)]B(x,y,z)}\bigl([x:A]\varphi(x),a,a,\rr(a)\bigr)\;=\;\varphi(a):B\bigl(a,a,\rr(a)\bigr)}$}
\end{prooftree}

\noindent Thus, using a trivial identity proof $\rr(a)$\/ to build a $J$-term does simply give back $\varphi(a)$.

To illustrate the use of the rules for derivations of judgements, we
give an example of a derivation which shows that the result of
applying a term to two propositionally equal terms results in propositionally equal terms.

\begin{example}\label{ex:prop_app}
Let $\tau \simeq \tau':\prod_{v:S}T(v)$, and let $\sigma:S$\/ be derivable. Then $\app(\tau,\sigma) \simeq \app(\tau',\sigma)$\/ is also derivable.
Indeed, consider the following derivation:

\vspace{.25cm}
\begin{prooftree}
  \AxiomC{$\judge{x,y:\prod_{v:S}T(v),z:\idn{\prod_{v:S}T(v)}{}(x,y)}{\idn{T(\sigma)}{}(\app(x,\sigma),\app(y,\sigma)):\type}$}
  \noLine
  \UnaryInfC{$\judge{x:\prod_{v:S}T(v)}{\rr(\app(x,\sigma)):\idn{T(\sigma)}{}(\app(x,\sigma),\app(x,\sigma))}$}
  \noLine
  \UnaryInfC{$\judge{}{f:\idn{\prod_{v:S}T(v)}{}(\tau,\tau')}$}
  \RightLabel{$\id{}$ elim.}
  \UnaryInfC{$\judge{}{\Jelim\bigl([x:\prod_{v:S}T(v)]\rr(\app(x,\sigma)),\tau,\tau',f):\idn{T(\sigma)}{}(\app(\tau,\sigma),\app(\tau',\sigma))}$}
\end{prooftree}

\vspace{.25cm}
\noindent Here, $f$\/ is a term witnessing the propositional identity $\tau \simeq \tau'$. Of course, the two other premises have to be derived as well,
but this is straightforward. 

Similarly we may derive from $\sigma \simeq \sigma'$\/ that $\app(\tau,\sigma) \simeq \app(\tau,\sigma')$.
\end{example}


\subsection{Natural numbers}\label{subsec:natural}

So far we have discussed only methods to construct new types and terms
from ones already present. It is common to introduce as a basic type the type
$\NN$\/ of natural numbers, and to add axioms
\begin{prooftree}
  \AxiomC{}
  \RightLabel{$\NN$ introduction (i)}
  \UnaryInfC{$\zero:\NN$}
\end{prooftree}
\vspace{.25cm}
\begin{prooftree}
  \AxiomC{$n:\NN$}
  \RightLabel{$\NN$ introduction (ii)}
  \UnaryInfC{$\successor(n):\NN$}
\end{prooftree}
which allow us to construct the standard numerals. Since the type of
natural numbers will not play a central role in this paper we refer to
the appendix for the precise
formulation of the elimination rule (expressing the possibility of
defining terms by recursion) and the conversion rules.
We do point out however that aside from the standard numerals the
theory may prove the existence of other, non-standard, numerals as
well.

\subsection{Theories and extensions}\label{subsec:theories}

We shall denote by $\mathbb{T}_\omega$\/ the system having all of
the above constructors and rules, including those for the type of natural numbers (for a complete description see the appendix). 
By a \myemph{type theory} we shall mean any
extension of the basic system $\mathbb{T}_\omega$\/ obtained by adding axioms and possibly also inference rules. The axioms are judgements
which may assert the existence of basic types or terms, or may assert the equality between certain types or terms. Possible additional inference rules include
the so-called \myemph{truncation-} and \myemph{reflection rules}, which express triviality of certain identity types. See Section~\ref{sec:monad} for a discussion of these rules.

Given two type theories $\mathbb{T}$\/ and $\mathbb{T}'$, we say that $\mathbb{T}'$\/ is an \myemph{extension} of $\mathbb{T}$\/ when every judgement which is
derivable in $\mathbb{T}$\/ is also derivable in $\mathbb{T}'$. Notation: $\mathbb{T} \subseteq \mathbb{T}'$. Thus by our definitions, $\mathbb{T}_\omega$\/ is
the smallest type theory.  

\subsection{Expressions}\label{subsec:raw}

Because the types and terms of such theories are defined
simultaneously, in order to formally specify the syntax of the theory
it is convenient to first define inductively a class of
\emph{expressions} --- which need not satisfy any typing conventions
--- from which the genuine syntactical data of the theory is then
extracted via the rules given above (and stated in full in Appendix
\ref{sec:types_ap}).  For example, in order to formally define the theory
$\mathbb{T}_{\omega}$ we first fix a countable set $V$ of (untyped)
variables and then define the class of \myemph{expressions of $\mathbb{T}_{\omega}$}, denoted
$\myexp(\mathbb{T}_{\omega})$, by
\begin{itemize}
\item $v$ is in $\myexp(\mathbb{T}_{\omega})$, for any $v$ in $V$;
\item $0$ and $\NN$ are in $\myexp(\mathbb{T}_{\omega})$;
\item
  $\successor(p),\app(p,q),\pair(p,q),\rr(p),\rec(p,q,r),\RSigma(p,q),\Jelim(p,q,r,s)$
  and $\lambda_{p:q}.r$ are in $\myexp(\mathbb{T}_{\omega})$ when $p,q,r$
  and $s$ are;
\item $\idn{p}{}(q,r),\prod_{p:q}r$ and $\sum_{p:q}r$ are in
  $\myexp(\mathbb{T}_{\omega})$ when $p,q$ and $r$ are.
\end{itemize}
Thus, the expressions are generated by applying all term- and type
constructors without regard for well-typedness. The derivation rules
of the type theory may then be regarded as carving out from this set
of all expressions those which are well-formed and well-typed.
The syntax of other the theories extending
$\mathbb{T}_{\omega}$ that we consider later is similarly specified in
this way with the evident modifications to the definition of the
expressions.  Moreover, because the expressions are inductively
generated it follows that the sets of the form $\myexp(-)$
possess an obvious universal property.

\subsection{Context morphisms}\label{subsec:context_morphisms}

Recall that if $\Gamma$ and 
\begin{align*}
  \Delta & = \bigl(x_{1}:A_{1},\ldots,x_{n}:A_{n}(x_{1},\ldots,x_{n-1})\bigr)
\end{align*}
are contexts, then a \myemph{context morphism} $a:\Gamma\to\Delta$ is
a sequence of terms
\begin{align*}
  \judge{\Gamma}{a_{1}:A_{1}},\judge{\Gamma}{a_{2}:A_{2}(a_{1})},\ldots,\judge{\Gamma}{a_{n}:A_{n}(a_{1},\ldots,a_{n-1})}.
\end{align*}
There is a category of contexts with arrows the context morphisms
(cf.~\cite{Hofmann:SSDT}).

\subsection{Globular sets}\label{subsec:globular}

Globular sets are structures which form the basis for several definitions of higher dimensional category. One way to think of a globular set is
as a higher dimensional graph: not only are there vertices and edges between the vertices, but one has edges between edges, and so on.
Formally, a globular set $G$\/ is a tuple $(G_n, s_n, t_n)_{n \in \mathbb{N}}$, where each $G_n$\/ is a set, and where $s_n, t_n:G_{n+1} \to G_n$\/ are functions
subject to the \myemph{globular identities}
\begin{equation} \label{eqn:glob}
d_nd_{n+1}=d_ns_{n+1}; \qquad s_ns_{n+1}=s_nd_{n+1} 
\end{equation}
for $d=s,t$.  Elements of $G_n$\/ are referred to as \myemph{n-cells}, and are said to have \myemph{dimension n}. The maps $s_n$\/ and $t_n$\/ are called \myemph{source} and
\myemph{target} maps, respectively.  

If $G$\/ is a globular set for which $G_n=\emptyset$\/ for all $n>1$, then we may simply regard $G$\/ as a (directed) graph. If there exist elements of higher
dimension, then the globular identities ensure that the source $s_n(x)$\/ and target $t_n(x)$\/ for such an $n$-dimensional edge are a parallel pair of edges
of dimension $n-1$.

Because it is often convenient,
given a $n$-cell $\alpha$ of a globular set $G$, to be able to refer
to the result of iteratively taking the source or target of $\alpha$
we introduce the notation $\alpha^{j}_{0},\alpha^{j}_{1}$ for these
corresponding $j$-cells.  Explicitly, for
$0\leq j\leq n-1$,
\begin{align*}
  \alpha^{j}_{i} & := \begin{cases}
    s_{j}\circ \cdots\circ s_{(n-1)}(\alpha) &
    \text{ if }i=0\\
    t_{j}\circ\cdots\circ t_{(n-1)}(\alpha) &
    \text{ if }i=1.
    \end{cases}
\end{align*}
By the globular identities, $\alpha^{j}_{0}$ and $\alpha^{j}_{1}$ are the
only elements of $G_{j}$ which are obtainable from $\alpha$ by applying the source and target maps.

A globular set $G$\/ is said to be \myemph{reflexive} if it comes equipped with a family of maps $i_n:G_n \to G_{n+1}$, such that 
\begin{equation} \label{eqn:rglob}
t_ni_n=1=s_ni_n
\end{equation} 
We think of $i_n(x)$\/ as the identity edge from $x$\/ to itself. In this paper we shall be working with reflexive globular sets only.
For readability we often omit the dimension from the source, target and identity maps of a globular set.

A \myemph{morphism of globular sets} $f:G \to H$\/ is a family of functions $f_n:G_n \to H_n$\/ which commute with the source and target maps.  Globular sets and their morphisms form a category denoted $\gsets$. For reflexive globular sets we also require that the $f_n$\/ commute with the identity maps; this gives
a category $\rgsets$.

There is a functor $\Delta:\sets \to \rgsets$\/ which takes a set $A$\/ to the constant globular set with $\Delta(A)_n=A$. A globular set which is isomorphic
to one of the form $\Delta(A)$\/ will be called \myemph{constant}. The functor $\Delta$\/ has a left adjoint
$\pi_0:\rgsets \to \sets$; this functor assigns to a globular set $G$\/ its set of connected components 
\[ \pi_0(G)=G_0/\!\!\sim\]
where the equivalence relation on 0-dimensional elements is generated by
\[ x \sim y \; \Leftrightarrow \exists f \in G_1. s(f)=x, t(f)=y. \]
We may express this as a (reflexive) coequalizer diagram:
\[ 
\xymatrix{
G_1 \ar@<.5ex>[r]^s \ar@<-.5ex>[r]_t& G_0 \ar@{->>}[r] & \pi_0(G). }
\]
The composite $\Delta\pi_0:\rgsets \to \rgsets$\/ is an idempotent monad, to which we refer as \myemph{0-truncation}. Often it will be convenient to 
identify the essential image of this functor (the constant globular sets) with the category of sets.

We may also truncate a globular set at dimension 1: in this case we replace the category of sets by the category $\rgrph$\/ of directed reflexive graphs.
There is a functor $\rgsets \to \rgrph$\/ which assigns to a globular set $G$\/ the graph whose vertex set is $G_0$\/ and whose edge set is $G_1/\sim$,
where two edges $h,k$\/ satisfy $h \sim k$\/ if there is an $\alpha \in G_2$\/ with $s(\alpha)=h, t(\alpha)=k$.

In the other direction, any directed reflexive graph $G$\/ gives a globular set with which is the same as $G$\/ in dimensions 0 and 1, and is degenerate above 
dimension 1. The composite functor $\rgsets \to \rgrph \to \rgsets$\/ will be called \myemph{1-truncation}, and a globular set in the essential image
of this functor will be said to be \myemph{1-truncated}. We shall
often identify the subcategory of 1-truncated globular sets with the
category of graphs.


\section{The Martin-L\"{o}f complex monad}\label{sec:monad}

The goal of this section is to state the formal definition of
Martin-L\"{o}f complexes.  Because Martin-L\"{o}f complexes are defined
to be algebras for a monad on the category of reflexive globular sets the
principal matter addressed here is the construction of the appropriate
monad.  The monad for the theory $\mathbb{T}_{\omega}$ obtained by the
construction below essentially corresponds to the monad obtained from the operadic
constructions due to van den Berg and Garner \cite{Berg:TWOG} and
Lumsdaine \cite{Lumsdaine:WOITT}, who show that the algebras are weak
omega-groupoids.  It is worth emphasizing that, because the converse
seems not to hold, the problem of determining precisely the higher-dimensional
structure of the algebras for these monads remains open.  It is to the
solution of this problem that the results of the present paper contribute.

Because we will be interested in algebras for the monad generated
by theories, such as the theories $\mathbb{T}_{n}$ described in Section
\ref{subsec:mlcxs} below, which extend
$\mathbb{T}_{\omega}$ the description of the monad involved in the
definition of Martin-L\"{o}f complexes will be described for an
arbitrary extension of $\mathbb{T}_{\omega}$.  As such, throughout
this section $\mathbb{T}$ is assumed to be an arbitrary theory
extending $\mathbb{T}_{\omega}$.  Finally, observe that although we choose to
work with reflexive globular sets, the construction of the
monad can be modified to yield a corresponding monad for globular
sets.

\subsection{Notation for iterated identity types and other conventions}

In order to most efficiently (and readably) state some of the
additional principles for identity types that we consider it
is useful to introduce notation for iterated identity types.
Fixing a type $A$ together with terms $a,b:A$ in some ambient context, we
introduce the (at this stage superfluous) notation
\begin{align*}
  \idn{A}{0} & := A,\text{ and}\\
  \idn{A}{1}(a,b) & := \idn{A}{}(a,b).
\end{align*}
In general, assuming given terms
\begin{align*}
  \judge{}{a_{n+1},b_{n+1}:\idn{A}{n}(a_{1},b_{1};\cdots;a_{n},b_{n})},
\end{align*}
we define
\begin{align*}
  \idn{A}{n+1}(a_{1},b_{1};\cdots;a_{n},b_{n};a_{n+1},b_{n+1}) & := \idn{\idn{A}{n}(a_{1},b_{1};\cdots;a_{n},b_{n})}{}(a_{n+1},b_{n+1}).
\end{align*}
In the sequel we will be dealing extensively with sets of terms from
various theories extending $\mathbb{T}_{\omega}$.  We adopt the
convention that such terms are always assumed to be identified modulo
definitional equality and $\alpha$-equivalence.

As a notational convenience we adopt the convention of, given
a reflexive globular set $G=(G_{n})_{n\geq 0}$, writing $G$ for the set
$\sum_{n\geq 0}G_{n}$.

\subsection{The reflexive globular set generated by a type}\label{subsec:glob_to_type}

Fix a type $A$ in $\mathbb{T}$.  It is possible that $A$ is a type in context, yet we will
assume that $A$ is a type in the empty context.  The case where the
context is non-empty is obtained in essentially the same way, and so
this is a reasonable simplification.  We will now construct a
reflexive globular set denoted by $\globfor{A}_{\mathbb{T}}$ and called the
\myemph{reflexive globular set generated by $A$ (in $\mathbb{T}$)}.
When the theory $\mathbb{T}$ is fixed we will omit the subscript and
write simply $\globfor{A}$.  This
construction will be carried out in such a way that the following
conditions are satisfied:
\begin{enumerate}
  \item Each element of $\globfor{A}_{n}$ is a tuple of $(2n+1)$ elements
    of the set of terms of $\mathbb{T}$.
  \item If both $(\vec{\alpha};\beta)$ and $(\vec{\alpha};\beta')$
    are in $\globfor{A}_{n}$, then
    $\judge{}{\idn{A}{n+1}(\vec{\alpha};\beta,\beta'):\type}$ is  derivable
    in $\mathbb{T}$.
  \item The source and target maps $s,t:\globfor{A}_{n+1}\rightarrow\globfor{A}_{n}$
    must send a tuple $(\alpha_{0},\ldots,\alpha_{2n})$ to
    $(\alpha_{0},\ldots,\alpha_{2n-2})$ and
    $(\alpha_{0},\ldots,\alpha_{2n-3},\alpha_{2n-1})$,
    respectively.
\end{enumerate}
We begin by defining 
\begin{align*}
  \globfor{A}_{0} & := \{a\;|\;\judge{}{a:A}\},\\
  \globfor{A}_{1} & :=
  \bigl\{(a_{0},a_{1};\;\alpha)\;|\;a_{0},a_{1}\in\globfor{A}_{0}\text{ and }\judge{}{\alpha:A(a_{0},a_{1})}\bigr\},
\end{align*}
and the maps $s,t:\globfor{A}_{1}\rightarrow\globfor{A}_{0}$ are simply the projections
$\pi_{0},\pi_{1}$ sending $(a_{0},a_{1};\;\alpha)$ to $a_{0}$ and
$a_{1}$, respectively.  Assuming $\globfor{A}$ has been constructed up
to stage $n$, we define $\globfor{A}_{n+1}$ to be the following set
\begin{align*}
  \bigl\{(\vec{\alpha};\;\beta_{0},\beta_{1};\;\gamma)\;|\;(\vec{\alpha};\;\beta_{i})\in\globfor{A}_{n}\text{
    for }\;i=0,1,\text{
     and }\judge{}{\gamma:\idn{A}{n+1}(\vec{\alpha};\;\beta_{0},\beta_{1})}\bigr\}.
\end{align*}
The source and target maps $s,t:\globfor{A}_{n+1}\rightarrow\globfor{A}_{n}$ are given
by the projections 
\begin{align*}
  \xy
  {\ar@{|->}(0,0)*+{(\vec{\alpha};\;\beta_{0},\beta_{1};\;\gamma)};(25,0)*+{(\vec{\alpha};\;\beta_{i}),}};
  \endxy
\end{align*}
for $i=0$ and $i=1$, respectively.
\begin{lemma}
  \label{lemma:type_gives_gset}
  Given an extension $\mathbb{T}$ of $\mathbb{T}_{\omega}$ and a
  (closed) type $A$ of $\mathbb{T}$, the graded set $\globfor{A}$ described
  above is a reflexive globular set.
  \begin{proof}
    The maps $i:\globfor{A}_{n}\to\globfor{A}_{n+1}$ are obtained
    using reflexivity terms.  The equations for reflexive globular
    sets are then readily verified.
  \end{proof}
\end{lemma}

\subsection{The type theory associated to a reflexive globular set}\label{subsec:type_to_glob}

Not only does every type $A$ give rise to a reflexive globular set, but also
every reflexive globular set $G$ gives rise to a type theory $\mathbb{T}[G]$.
\begin{definition}
  \label{def:tt_gen_by_glob}
  Given a reflexive globular set $G$, the type theory $\mathbb{T}[G]$ \myemph{generated by
    $G$} (or \myemph{$\mathbb{T}$ with $G$ adjoined}) is obtained by
  augmenting $\mathbb{T}$ with the following additional symbols and
  rules:
  \begin{itemize}
  \item A basic type $\judge{}{\namefor{G}}$;
  \item Basic terms
    $\judge{}{\namefor{g}:\namefor{G}}$, for each vertex $g\in G_{0}$;
  \item Basic terms
    $\judge{}{\namefor{f}:\namefor{G}(\namefor{g},\namefor{h})}$, for
    each element $f\in G_{1}$ with $s(f)=g$ and $t(f)=h$;
  \item Basic terms
    \begin{align}
      \label{eq:new_id_judgements}
      \judge{}{\namefor{\alpha}:\namefor{G}^{\;n}\bigl(\namefor{\alpha^{0}_{0}},\namefor{\alpha^{0}_{1}};\;\namefor{\alpha^{1}_{0}},\namefor{\alpha^{1}_{1}};\;\cdots;\;\namefor{\alpha^{n-1}_{0}},\namefor{\alpha^{n-1}_{1}}\bigr)}
    \end{align}
    where $\alpha^{j}_{i}$ for $i=0,1$ and $0\leq j\leq n-1$ are as
    defined in Section \ref{subsec:globular}, for each $\alpha\in G_{n}$;
  \item New conversion rules:
    \begin{align*}
      \namefor{i(\alpha)}\;=\;\rr\bigl(\namefor{\alpha}\bigr):\namefor{G}^{n+1}\bigl(\ldots;\namefor{\alpha},\namefor{\alpha}\bigr)
  \end{align*}
  for every $\alpha\in G_{n}$.
  \end{itemize}  
\end{definition}
\begin{remark}
  As a matter of
  notation, we write $\judge[G]{\Gamma}{\mathcal{J}}$ to indicate that
  the judgement \mbox{$\judge{\Gamma}{\mathcal{J}}$} is derivable in
  $\mathbb{T}[G]$.  Finally, we also write
  $\myexp_{G}$ instead of the more cumbersome
  $\myexp(\mathbb{T}[G])$.  Also, when no confusion
  will result, we identify the symbol $\namefor{\tau}$ with $\tau$
  itself.  E.g., we write $f:\idn{G}{}(g,h)$ instead of
  the more cumbersome $\namefor{f}:\namefor{G}(\namefor{g},\namefor{h})$.
\end{remark}  
In subsequent sections it will be convenient to have at our disposal
techniques for constructing maps between the sets of expressions of
one type theory $\mathbb{T}[G]$ and another $\mathbb{T}[H]$, for $G$
and $H$ globular sets.  Along these lines, we make the following
observation.
\begin{lemma}
  \label{lemma:ext}
  Given globular sets $G$ and $H$, any function
  \begin{align*}
    \xy
    {\ar^{\varphi}(0,0)*+{G};(20,0)*+{\myexp_{H}}};
    \endxy
  \end{align*}
  has a unique extension
  $\hat{\varphi}:\myexp_{G}\to\myexp_{H}$,
  commuting with the operations from which the expressions
  are formed, such that the following diagram of sets commutes:
  \begin{align*}
    \xy
    {\ar@{..>}^{\hat{\varphi}}(0,15)*+{\myexp_{G}};(30,15)*+{\myexp_{H}}};
    {\ar^{i_{G}}(15,0)*+{G};(0,15)*+{\myexp_{G}}};
    {\ar_{\varphi}(15,0)*+{G};(30,15)*+{\myexp_{H}}};
    \endxy
  \end{align*}
  where $i_{G}$ is the map sending $g\in G_{n}$ to $\namefor{g}$.
\end{lemma}
Note that the basic type $\namefor{G}$ is sent by the extension
$\hat{\varphi}$ to $\namefor{H}$.  Of course, depending on the nature
of $\varphi$ the extension $\hat{\varphi}$ may or may not preserve
derivable judgements.  Such a $\hat{\varphi}$ will, however, commute with
substitution.  I.e., if $e(x)$ is an expression of $\mathbb{T}[G]$
with $x$ free, then, for any other expression $f$, 
\begin{align}
  \label{eq:ext_sub}
  \hat{\varphi}(e)[\hat{\varphi}(f)/x] & = \hat{\varphi}(e[f/x]).
\end{align}

\subsection{The induced monad on globular sets}\label{section:induced_monad}

We will now see that composing the foregoing processes
\begin{align*}
  G & \longmapsto\mathbb{T}[G],\text{ and}\\
  A:\type & \longmapsto\globfor{A},
\end{align*}
yields a monad $T$ on the category $\rglob$ of reflexive globular sets.  Given a globular set $G$,
\begin{align}
  \label{eq:def_type_monad}
  T(G) & := \globfor{\namefor{G}}.
\end{align}
Suppose given a map $\varphi:G\to H$ of globular sets.  
To see that this assignment is in fact functorial we begin by noting
that, by Lemma \ref{lemma:ext}, the map
\begin{align*}
  g & \longmapsto\namefor{\varphi(g)}
\end{align*}
for $g\in G_{n}$, possesses a canonical extension 
$\varphi_{*}:\myexp_{G}\to\myexp_{H}$.  I.e., in the notation of
Lemma \ref{lemma:ext},
\begin{align*}
  \varphi_{*} & :=\widehat{i_{H}\circ\varphi}.
\end{align*}
In order to be able to use $\varphi_{*}$ to define the action of $T$
on arrows we must first verify that it preserves derivable
judgements, where the action of $\varphi_{*}$ extends to judgements
in the obvious manner.
\begin{lemma}
  \label{lem:hat_preserves_judge}
  Suppose $\mathcal{J}$ is a judgement derivable in
  $\mathbb{T}[G]$, then $\varphi_{*}(\mathcal{J})$ is derivable in $\mathbb{T}[H]$.
  \begin{proof}
    The proof is a straightforward induction on the structure of derivations
    $\judge[G]{}{\mathcal{J}}$.  For example, suppose $\mathcal{J}$ is
    the conclusion $\judge[G]{\Gamma}{\lambda_{x:A}.b(x):\prod_{x:A}.B(x)}$
    of the introduction rule for dependent products.  Then we have by
    the induction hypothesis that
    \begin{align*}
      \judge[H]{\varphi_{*}(\Gamma),x:\varphi_{*}(A)}{\varphi_{*}(b(x)):\varphi_{*}(B(x))}.
    \end{align*}
    Applying the introduction rule in $\mathbb{T}[H]$ yields
    the appropriate judgement since
    \begin{align*}
      \varphi_{*}\bigl(\prod_{x:A}.B(x)\bigr) & = \prod_{x:\varphi_{*}(A)}.\varphi_{*}(B)(x),
    \end{align*}
    by definition of $\varphi_{*}$.  The only case which merits
    special attention are those judgements of the form (\ref{eq:new_id_judgements})
    which occur as axioms of $\mathbb{T}[G]$.  Such judgements are
    preserved by the fact that $\varphi$ is a map of globular sets.
  \end{proof}
\end{lemma}
\begin{lemma}
  \label{lem:monad_funct}
  The assignment (\ref{eq:def_type_monad}) is functorial
  $T:\rgsets\to\rgsets$.
  \begin{proof}
    Let
    \begin{align*}
      T(\varphi)(\alpha_{0},\alpha_{1},\;\cdots,\;\alpha_{2n})
      & := \bigl(\varphi_{*}(\alpha_{0}),\varphi_{*}(\alpha_{1}),\;\cdots,\;\varphi_{*}(\alpha_{2n})\bigr),
    \end{align*}
    for $\vec{\alpha}$ in $T(G)_{n+1}$.  That this definition makes sense
    follows from Lemma \ref{lem:hat_preserves_judge} and the definition of
    $\varphi_{*}$.  Trivially, $T(1_{G})=1_{T(G)}$.  To see that $T$
    is well behaved with respect to composition it suffices to show
    that, when given $\psi:H\to I$, we have
    \begin{align*}
      (\psi\circ\varphi)_{*} & =\psi_{*}\circ\varphi_{*}.
    \end{align*}
    For this we observe that on the generators $g\in G_{n}$,
    \begin{align*}
      \psi_{*}\bigl(\varphi_{*}(\namefor{g})\bigr)\; = \;\psi_{*}\bigl(\namefor{\varphi(g)}\bigr)\;
       =\;\namefor{\psi\circ\varphi(g)}\;
       = \;(\psi\circ\varphi)_{*}(\namefor{g}).
    \end{align*}
  \end{proof}
\end{lemma}
As a notational convenience we will often write elements
$\vec{\alpha}\in T(G)$ in terms of their boundaries.  I.e., we write
$\vec{\alpha}=(\alpha^{0}_{0},\alpha^{0}_{1};\ldots;\alpha^{n-1}_{0},\alpha^{n-1}_{1};\alpha)$
instead of $(\alpha_{0},\alpha_{1},\ldots,\alpha_{2n})$.
\begin{proposition}
  \label{prop:M_is_monad}
  $T:\rgsets\to\rgsets$ is the functor part of a
  monad.
  \begin{proof}
    Given a globular set $G$, the unit $\eta_{G}:G\to T(G)$ is the
    ``insertion of generators'' defined by setting
    \begin{align*}
      \eta_{G}(g) & := \bigl(\namefor{g^{0}_{0}},\namefor{g^{0}_{1}},\cdots,\namefor{g}\bigr)
    \end{align*}
    for $g\in G_{n}$.  This is a globular map which is natural in $G$
    by definition.  

    For the multiplication $\mu_{G}:T^{2}G\to TG$ we begin by
    defining
    $\tau_{G}:\myexp_{TG}\to\myexp_{G}$ to be
    the canonical extension, which exists by Lemma
    \ref{lemma:ext}, of the assignment $\pi:TG\to\myexp_{G}$ given by
    \begin{align*}
      (\alpha^{0}_{0},\alpha^{0}_{1};\ldots;\alpha^{n-1}_{0},\alpha^{n-1}_{1};\alpha)&\longmapsto\alpha,
    \end{align*}
    where $\vec{\alpha}$ is in $(TG)_{n}$.  That is, $\tau_{G}=\hat{\pi}$ is
    the canonical extension such that
    \begin{align*}
      \xy
      {\ar^{\tau_{G}}(0,15)*+{\myexp_{TG}};(30,15)*+{\myexp_{G}}};
      {\ar^{i_{TG}}(15,0)*+{TG};(0,15)*+{\myexp_{TG}}};
      {\ar_{\pi}(15,0)*+{TG};(30,15)*+{\myexp_{G}}};
      \endxy
    \end{align*}
    commutes.   As such, given $\vec{\alpha}$ in $(TG)_{n}$ as above,
    \begin{align*}
      \tau_{G}\bigl(\namefor{\vec{\alpha}}\bigr) & = \alpha,
    \end{align*}
    where this definition makes sense because $\alpha$ is
    itself a term of $\mathbb{T}[G]$.  We would like to show that
    $\tau_{G}$ preserves derivable judgements.  As in the proof of Lemma
    \ref{lem:hat_preserves_judge}, the non-trivial step is to verify
    that the axioms added in the formation of $\mathbb{T}[TG]$ are
    preserved.  That is, where $\vec{\alpha}$ is as above, we need to
    show that
    \begin{align*}
      \judge[TG]{}{\namefor{\vec{\alpha}}:\namefor{TG}^{n}\bigl(\namefor{(\vec{\alpha})^{0}_{0}},\ldots,\namefor{(\vec{\alpha})^{n-1}_{1}}\bigr)}
    \end{align*}
    implies the corresponding judgement in $\mathbb{T}[G]$.  But,
    we have that 
    \begin{align}\label{eq:tau_G}
      \tau_{G}\bigl(\namefor{(\vec{\alpha})^{i}_{j}}\bigr) & = \alpha^{i}_{j}.
    \end{align}
    Thus, we must show that
    \begin{align*}
      \judge[G]{}{\alpha:\namefor{G}^{n}(\alpha_{0}^{0},\ldots,\alpha^{n-1}_{1})}.
    \end{align*}
    However, this is a trivial consequence of the fact that
    $\vec{\alpha}$ is an element of $(TG)_{n}$.  Therefore $\tau_{G}$ preserves derivable
    judgements and we may define
    \begin{align*}
      \mu_{G}(\beta^{0}_{0},\ldots,\beta^{n-1}_{1},\beta) & := \bigl(\tau_{G}(\beta^{0}_{0}),\ldots,\tau_{G}(\beta^{n-1}_{1}),\tau_{G}(\beta)\bigr),
    \end{align*}
    for $\vec{\beta}$ in $(T^{2}G)_{n}$.  Since $\tau_{G}$ preserves
    valid judgements this gives a globular map which is natural in $G$.

    To see that the first unit law for monads is satisfied, let
    $\vec{\alpha}$ in $(TG)_{n}$ be given as above.  Then 
    \begin{align*}
      \mu_{G}\circ\eta_{TG}(\vec{\alpha}) & =
      \mu_{G}\bigl(\namefor{(\vec{\alpha})_{0}^{0}},\ldots,\namefor{(\vec{\alpha})^{n-1}_{1}},\namefor{\vec{\alpha}}\bigr)\\
      & =
      \biggl(\tau_{G}\bigl(\namefor{(\vec{\alpha})_{0}^{0}}\bigr),\ldots,\tau_{G}\bigl(\namefor{(\vec{\alpha})^{n-1}_{1}}\bigr),\tau_{G}\bigl(\namefor{\vec{\alpha}}\bigr)\biggr)\\
      & = \vec{\alpha},
    \end{align*}
    where the final equation is by (\ref{eq:tau_G}).  For the other
    unit law, observe that the following diagram commutes:
    \begin{align*}
      \xy
      {\ar^{i_{G}}(15,0)*+{G};(0,15)*+{\myexp_{G}}};
      {\ar_{\eta_{G}}(15,0)*+{G};(45,0)*+{TG}};
      {\ar^{(\eta_{G})_{*}}(0,15)*+{\myexp_{G}};(30,15)*+{\myexp_{TG}}};
      {\ar_{i_{TG}}(45,0)*+{TG};(30,15)*+{\myexp_{TG}}};
      {\ar^{\tau_{G}}(30,15)*+{\myexp_{TG}};(60,15)*+{\myexp_{G}}};
      {\ar_{\pi}(45,0)*+{TG};(60,15)*+{\myexp_{G}}};
      \endxy
    \end{align*}
    Thus, $\tau_{G}\circ(\eta_{G})_{*}\circ i_{G}=i_{G}$ and, by
    Lemma \ref{lemma:ext}, $\mu_{G}\circ T(\eta_{G})=1_{TG}$.

    Next, to see that the multiplication law is satisfied it suffices
    to prove that 
    \begin{align}
      \label{eq:glob_monad_goal}
      \tau_{G}\circ\tau_{TG} & = \tau_{G}\circ(\mu_{G})_{*}.
    \end{align}
    Given $\vec{\beta}=(\beta^{0}_{0},\beta^{0}_{1};\ldots;\beta^{n-1}_{0},\beta^{n-1}_{1};\beta)$ in $(T^{2}G)_{n}$ we have
    \begin{align*}
      \tau_{G}\circ\tau_{TG}\bigl(\namefor{\vec{\beta}}\bigr) & =
      \tau_{G}(\beta)\\
      & =
      \tau_{G}\biggl(\namefor{\bigl(\tau_{G}(\beta^{0}_{0}),\tau_{G}(\beta^{0}_{1}),\ldots,\tau_{G}(\beta)\bigr)}\biggr)\\
      & =
      \tau_{G}\bigl(\namefor{\mu_{G}(\vec{\beta})}\bigr)\\
      & = \tau_{G}\circ(\mu_{G})_{*}\bigl(\namefor{\vec{\beta}}\bigr).
    \end{align*}
    Thus, by Lemma \ref{lemma:ext},
    (\ref{eq:glob_monad_goal}) holds.
  \end{proof}
\end{proposition}
\begin{example}
  Suppose $g$ is a vertex of $G$, then $\namefor{(\namefor{g})}$ is
  likewise a vertex of $T^{2}G$.  The multiplication $\mu_{G}$ acts on
  such a vertex by removing the outermost $\namefor{-}$.  I.e., 
  \begin{align*}
    \mu_{G}\bigl(\namefor{(\namefor{g})}\bigr) & = \namefor{g}.
  \end{align*}
  Similarly, if $f$ is in $G_{n}$, then
  \begin{align*}
    \mu_{G}\bigl(\namefor{(\namefor{f^{0}_{0}},\namefor{f^{0}_{1}},\ldots,\namefor{f^{n-1}_{1}},\namefor{f})}\bigr)
    & = (\namefor{f^{0}_{0}},\ldots,\namefor{f}).
  \end{align*}
  The action of $\mu_{G}$ on composite terms (constructed out of the
  basic terms of $\mathbb{T}[TG]$ using the rules of $\mathbb{T}$) is
  then to go through the term recursively removing occurrences of
  $\namefor{-}$.  Thus, the unit acts by adding $\namefor{-}$ and the
  multiplication acts by removing it.
\end{example}

\subsection{Martin-L\"{o}f complexes and other categories of algebras}\label{subsec:mlcxs}

It is possible to extend Proposition \ref{prop:M_is_monad} by  
allowing the extension $\mathbb{T}$ of $\mathbb{T}_{\omega}$
employed in the construction to vary.  We denote by $\extomega$ the
category of all extensions of $\mathbb{T}_{\omega}$.  I.e., the
objects of $\extomega$ are dependent type theories extending
$\mathbb{T}_{\omega}$ (where we only allow those extensions
obtained by the addition of set-many new symbols and rules).  A
morphism $\mathbb{T}\to\mathbb{T}'$ in $\extomega$ is an inclusion of
theories (i.e., such a morphism exists whenever $\mathbb{T}'$ extends
$\mathbb{T}$).  We also denote by $\monadsin(\rglob)$ the category of
monads on $\rglob$ (regarded as monoids in $[\rglob,\rglob]$).
\begin{lemma}
  \label{lem:extension_functor}
  The construction of a monad on $\rglob$ from an extension of
  $\mathbb{T}_{\omega}$ from Section \ref{section:induced_monad} gives
  the action on objects of a functor
  \begin{align*}
   \mathcal{T}:\extomega\longrightarrow\monadsin(\rglob).
  \end{align*}
  \begin{proof}
    Assume given theories $\mathbb{T}$ and $\mathbb{T}'$ in
    $\extomega$ such that $\mathbb{T}'$ is an extension of
    $\mathbb{T}$.  We will now describe the induced natural transformation
    $\xi:T\to T'$, where we write $T$ and $T'$ as abbreviations for
    $\mathcal{T}(\mathbb{T})$ and $\mathcal{T}(\mathbb{T}')$,
    respectively.  Given a
    reflexive globular set $G$ and an element $\vec{\alpha}=(\alpha^{0}_{0},\ldots,\alpha^{n-1}_{1},\alpha)$ of
    $T(G)_{n}$, we note that since $\mathbb{T}'$ extends $\mathbb{T}$
    it follows that each component of the list $\vec{\alpha}$ is also
    a term of $\mathbb{T}'[G]$.  Moreover, all of these terms
    necessarily possess the appropriate boundaries so that
    $\vec{\alpha}$ is also an element of $T'(G)_{n}$.  As such, we may
    simply define $(\xi_{G})_{n}:T(G)_{n}\to T'(G)_{n}$ to be the map
    which sends any $\vec{\alpha}$ as above to itself (now regarded as
    a list of terms from $\mathbb{T}'[G]$).  This clearly describes a
    map of reflexive globular sets which is clearly
    $\xi$ is natural and that it commutes with the multiplication and
    unit maps for $T$ and $T'$.  Finally, it is trivial to see that,
    with this definition $\mathcal{T}$ is functorial.
  \end{proof}
\end{lemma}
\noindent The specific extensions of $\mathbb{T}_{\omega}$ to which we would like
to apply Lemma \ref{lem:extension_functor} are obtained by
augmenting $\mathbb{T}_{\omega}$ by axioms that force the identity
types to be trivial once they have been iterated sufficiently many
times.  To begin with, recall that the \myemph{reflection rule} for identity types is the
principle which states that all identity types are
trivial in the sense that
\begin{prooftree}
  \AxiomC{$\judge{}{a,b:A}$}
  \AxiomC{$\judge{}{p:\idn{A}{}(a,b)}$}
  \RightLabel{Reflection}
  \BinaryInfC{$\judge{}{a=b:A}$}
\end{prooftree}
\noindent Higher-dimensional
generalizations of this rule are then given by ``truncating'' the
identity types only after they have been iterated a certain number of
times.  Explicitly, the \myemph{$n$-truncation rule} is stated as follows:
\begin{prooftree}
  \AxiomC{$\judge{}{a_{n+1},b_{n+1}:\idn{A}{n}(a_{1},b_{1};\cdots;a_{n},b_{n})}$}
  \AxiomC{$\judge{}{p:\idn{A}{n+1}(a_{1},b_{1};\cdots;a_{n+1},b_{n+1})}$}
  \RightLabel{TR$_{n}$}
  \BinaryInfC{$\judge{}{a_{n+1}=b_{n+1}:\idn{A}{n}(a_{1},b_{1};\cdots;a_{n},b_{n})}$}
\end{prooftree}
With these rules at our disposal we are able to describe the type
theories extending $\mathbb{T}_{\omega}$ with which we will be
concerned.  Explicitly, for $n\geq 0$, the theory $\mathbb{T}_{n}$ is
defined to be the result of adding to $\mathbb{T}_{\omega}$ the (instances of the)
principle $\text{TR}_{n}$.  These theories then arrange themselves
according to the following hierarchy of theories:
\begin{align}\label{eq:hierarchy}
  \mathbb{T}_{\omega}\;\subseteq\;\cdots\;\subseteq\;\mathbb{T}_{n+1}\;
  \subseteq\;\mathbb{T}_{n}\;\subseteq\;
  \cdots\;\subseteq\;\mathbb{T}_{1}\;\subseteq\;
  \mathbb{T}_{0},
\end{align}
since $\text{TR}_{m}$ clearly implies $\text{TR}_{n}$, when $m<n$.
The theory $\mathbb{T}_{0}$ is also known as \emph{extensional type
  theory} as contrasted with the \emph{intensional type theory}
$\mathbb{T}_{\omega}$.
\begin{definition}\label{def:MLCxs}
  Denote by $M_{\omega}$ the monad $\mathcal{T}(\mathbb{T}_{\omega})$.
  A reflexive globular set $G$ is a \myemph{Martin-L\"{o}f complex} (or 
  \myemph{ML-complex}) if it is an algebra for $M_{\omega}$.  We write $\mlcxs$
  for the Eilenberg-Moore category consisting of $M_{\omega}$-algebras and
  homomorphisms thereof.  Similarly, we denote by $\mlcxs_{n}$ the
  category of $M_{n}$-algebras for $n=0,1,2,\ldots$, where $M_{n}$ denotes
  the monad $\mathcal{T}(\mathbb{T}_{n})$.
\end{definition}
Corresponding to the hierarchy of theories (\ref{eq:hierarchy}) we obtain, by
Lemma \ref{lem:extension_functor}, the following sequence of inclusions of
categories:
\begin{align*}
  \mlcxs_{0}\longrightarrow\cdots\longrightarrow\mlcxs_{n}\longrightarrow\mlcxs_{n+1}\longrightarrow\cdots\longrightarrow\mlcxs_{\omega}
\end{align*}
and it is our goal to understand how these categories relate to the
hierarchy of categories of homotopy types discussed in Section
\ref{sec:intro}.  

\subsection{Connection between truncation and other rules}\label{subsec:truncation}

The truncation principles $\text{TR}_{n}$ are related to several other
type theoretic principles which we employ occasionally in the sequel.  For example consider the
following $n$-dimensional generalization of the principle of \myemph{(definitional) uniqueness of
  identity proofs}:
\begin{prooftree}
  \AxiomC{$\judge{}{a_{n+1},b_{n+1}:\idn{A}{n}(a_{1},b_{1};\cdots;a_{n},b_{n})}$}
  \RightLabel{UIP$_{n}$}
  \UnaryInfC{$\judge{}{a_{n+1}=b_{n+1}:\idn{A}{n}(a_{1},b_{1};\cdots;a_{n},b_{n})}$}  
\end{prooftree}
The question whether $\text{UIP}_{1}$ --- or the variant where the
definitional equality occurring in the conclusion is replaced by a
propositional equality --- is derivable in $\mathbb{T}_{\omega}$ was
one of the motivations for the original groupoid model due to Hofmann
and Streicher \cite{Hofmann:GITT}.  In particular, the groupoid model
shows that neither $\text{UIP}_{1}$ nor the propositional version are
derivable in $\mathbb{T}_{\omega}$.

Another related principle is the $n$-dimensional \myemph{ordinary unit principle}
\begin{prooftree}
  \AxiomC{$\judge{}{a_{n+1}:\idn{A}{n}(a_{1},b_{1};\cdots;a_{n},b_{n})}$}
  \AxiomC{$\judge{}{p:\idn{A}{n+1}(a_{1},b_{1};\cdots;a_{n+1},a_{n+1})}$}
  \RightLabel{OUP$_{n}$}
  \BinaryInfC{$\judge{}{p\;=\;\rr(a_{n+1}):\idn{A}{n+1}(a_{1},b_{1};\cdots;a_{n+1},a_{n+1})}$}
\end{prooftree}

\noindent Whereas the uniqueness of identity proofs principles can be thought of as
requiring that the identity types are preorders above a given
dimension, the ordinary unit rules indicate that all loops (above
certain dimensions) are necessarily identities.  

The truncation and ordinary unit principles have been considered
previously by Garner in \cite{Garner:TDMTT} and by Warren in
\cite{Warren:PhD}.  The relation between the truncation, uniqueness of
identity proofs and ordinary unit principles are clarified in the
following lemma (the idea for the proof of which comes essentially
from results, which are not ``stratified'' in the way considered here,
from \cite{Streicher:IIITT}).
\begin{lemma}
  \label{lem:truncation_UIP}
  Assuming the rules of $\mathbb{T}_{\omega}$ and the usual rules for
  identity types, the
  following implications hold:
  \begin{enumerate}
  \item $\textnormal{TR}_{n}$ implies $\textnormal{OUP}_{n}$.
  \item $\textnormal{TR}_{n}$ implies $\textnormal{UIP}_{n+1}$.
  \item $\textnormal{UIP}_{n}$ implies $\textnormal{TR}_{n}$.
  \end{enumerate}
  for $n\geq 0$.
  \begin{proof}
    For (1), let a term $a_{n+1}$ of type
    $\idn{A}{n}(a_{1},b_{1};\cdots;a_{n},b_{n})$ and a ``loop'' $p$ of
    type $\idn{A}{n+1}(a_{1},b_{1};\cdots;a_{n+1},a_{n+1})$ be given.
    Then, by $\textnormal{TR}_{n}$ it suffices to show that
    \begin{align*}
     \judge{}{p\;\homotopic\;\rr(a_{n+1})\;:\idn{A}{n+1}(a_{1},b_{1};\cdots;a_{n+1},a_{n+1})}.
    \end{align*}
    To this end, observe that, by $\textnormal{TR}_{n}$, the type
    \begin{align*}
      \judge{x,y:\idn{A}{n}(a_{n},b_{n}),z:\idn{A}{n+1}(x,y)}{\idn{A}{n+2}\bigl(z,\rr(x)\bigr)}
    \end{align*}
    is derivable.  Therefore, the elimination rule for identity types gives us
    \begin{align*}
      \Jelim(\rr(x),a_{n+1},a_{n+1},p):\idn{A}{n+2}\bigl(p,\rr(a_{n+1})\bigr),
    \end{align*}
    as required.  (This proof essentially shows that Streicher's $K$
    rule \cite{Streicher:IIITT} is derivable for identity types of the
    form $\idn{A}{n+1}$.)

    Suppose, for the proof of (2), that we are given terms $a_{n+2}$
    and $b_{n+2}$ of type
    $\idn{A}{n+1}(a_{1},b_{1};\cdots\;a_{n+1},b_{n+1})$.  Then, by
    $\text{TR}_{n}$, $a_{n+1}=b_{n+1}$.  By (1) it follows that
    $\text{OUP}_{n}$ holds and therefore we obtain
    \begin{align*}
      a_{n+2} \; =\; \rr(a_{n+1}) \; = \; b_{n+2},
    \end{align*}
    as required.

    Finally, (3) holds trivially.
  \end{proof}
\end{lemma}

\subsection{Skeletal terms}\label{section:skeleta}
In this section we introduce a convenient technical tool which will later facilitate reasoning
about the type theories at hand and their models.

Let us denote by $\mathbb{T}_{\kappa}^{+}$ the theory
$\mathbb{T}_{\kappa}[0]$ for $\kappa=0,1,\ldots,\omega$.  This theory
is the same as $\mathbb{T}_{\kappa}$ except that it also has a new
distinguished type symbol $\namefor{0}$.  All of the theories of the
form $\mathbb{T}_{\kappa}[G]$ are extensions of
$\mathbb{T}_{\kappa}^{+}$.
\begin{definition}
  A judgement $\mathcal{J}$ of $\mathbb{T}_{\kappa}[G]$ is
  \myemph{skeletal} if no basic term symbols (coming from $G$) occur
  in $\mathcal{J}$.
\end{definition}
If $\mathcal{J}$ is a skeletal judgement in $\mathbb{T}_{\kappa}[G]$,
then there is a corresponding judgement $\mathcal{J}^{+}$ in
$\mathbb{T}_{\kappa}^{+}$ obtained by replacing each occurrence of the
basic type $\namefor{G}$ by $\namefor{0}$.  Similarly, any judgement
$\mathcal{J}$ has a corresponding translation $\mathcal{J}^{G}$ into
any $\mathbb{T}_{\kappa}[G]$ by replacing each occurrence of
$\namefor{0}$ by $\namefor{G}$.  We have the following basic
observation about the derivability of skeletal judgements:
\begin{lemma}
  A skeletal judgement $\mathcal{J}$ is derivable in
  $\mathbb{T}_{\kappa}[G]$ if and only if $\mathcal{J}^{+}$ is
  derivable in $\mathbb{T}_{\kappa}^{+}$.
\end{lemma}
Obviously the analogous statement which says that $\mathcal{J}$ is
derivable in $\mathbb{T}_{\kappa}^{+}$ if and only if
$\mathcal{J}^{G}$ is derivable in $\mathbb{T}_{\kappa}[G]$ also
holds.  Accordingly, we will henceforth not distinguish between the
judgements $\mathcal{J}$, $\mathcal{J}^{+}$ and $\mathcal{J}^{G}$.

Assume that $H$ is a finite reflexive globular set and define a
context $\Delta_{H}$ in $\mathbb{T}_{\kappa}^{+}$ as follows.
$\Delta_{H}$ consists of:
\begin{itemize}
\item For each $0$-cell $a$ of $H$, there is a distinct variable
  declaration $v_{a}:\namefor{0}$.
\item For each $(n+1)$-cell $f$ of $H$, there is a distinct variable
  declaration $v_{f}:\idn{\namefor{0}}{n}(v_{s(f)},v_{t(f)})$.
\end{itemize}
Because $H$ is finite this determines a well-defined context.  

Fix some $\kappa=0,1,\ldots,\omega$.  We now define a new category $\CC^{\kappa}_{H}$ as follows.
\begin{definition}
  A \myemph{context
    relative to $H$} is a skeletal context $\Gamma$ extending
  $\Delta_{H}$.  Given contexts $\Gamma$ and $\Theta$ relative to $H$,
  a \myemph{context morphism $\sigma:\Gamma\to\Theta$ relative to $H$}
  is a skeletal context morphism such that 
  \begin{align*}
    \xy
    {\ar^{\sigma}(0,10)*+{\Gamma};(20,10)*+{\Theta}};
    {\ar(0,10)*+{\Gamma};(10,0)*+{\Delta_{H}}};
    {\ar(20,10)*+{\Theta};(10,0)*+{\Delta_{H}}};
    \endxy
  \end{align*}
  commutes.
\end{definition}
The category $\CC^{\kappa}_{H}$ has as objects contexts in $\mathbb{T}_{\kappa}^{+}$ relative to $H$ and as
arrows context morphisms in $\mathbb{T}_{\kappa}^{+}$ relative to $H$.  
In fact, we have
a comprehension category (cf.~\cite{Jacobs:CLTT}) with base $\CC^{\kappa}_{H}$ and with fibration
$P^{\kappa}_{H}:\TT^{\kappa}_{H}\to\CC^{\kappa}_{H}$ determined by letting the fiber
$\TT^{\kappa}_{H}(\Gamma)$ consists of the skeletal types in context $\Gamma$.
This determines a split Grothendieck fibration since skeletal types
are stable under skeletal substitutions and there is an obvious
comprehension map $\chi^{\kappa}_{H}:\TT^{\kappa}_{H}\to(\CC^{\kappa}_{H})^{\to}$ which sends a
skeletal type $\judge{\Gamma}{A}$ to the dependent projection
$(\Gamma,x:A)\to\Gamma$.  

\begin{lemma}
The comprehension category $\CC^{\kappa}_{H}$\/ is a model of $\mathbb{T}_{\kappa}[H]$.
\end{lemma}

\begin{proof}
There is an obvious forgetful functor $\CC^{\kappa}_{H} \to \CC_{\mathbb{T}_{\kappa}[H]}$, the usual
 syntactic model of $\mathbb{T}_{\kappa}[H]$,
and this functor preserves the comprehension category structure. 
All type and term formation operations respect skeletalness w.r.t. $\Delta_H$. 
Therefore, the category $\CC^{\kappa}_{H}$\/ also supports
dependent products, sums, natural numbers and identity types, and the forgetful functor creates these.

Now $\CC^{\kappa}_{H}$\/ becomes a model of $\mathbb{T}_{\kappa}[H]$\/
by interpreting a basic term $a$\/ as the
variable $v_{a}$.  
\end{proof}

Let us write $\textnormal{sk}_{H}(\mathcal{J})$ for the
interpretation of a judgement $\mathcal{J}$ of
$\mathbb{T}_{\kappa}[H]$.  It follows, by induction on derivations,
that the substitution $\sigma_{H}$ given by $v_{a}\mapsto\namefor{a}$
satisfies
\begin{align}\label{eq:skeleton}
  \textnormal{sk}_{H}(\mathcal{J})[\sigma_{H}] & \equiv \mathcal{J}.
\end{align}
Indeed, $\sigma_{H}$ induces a morphism of models
$\CC_{H}^{\kappa}\to\CC_{\mathbb{T}_{\kappa}[H]}$ where again
$\CC_{\mathbb{T}_{\kappa}[H]}$ is the syntactic model.  This morphism is in fact an 
isomorphism of models, with inverse induced by $\text{sk}_{H}$\/ (which sends contexts and
context morphisms to their ``skeletalized" counterparts, i.e. it replaces each basic term by the
appropriate variable). 

In the case where $H$\/ is no longer finite matters become less straightforward.
\begin{lemma}
  \label{lemma:skeleta}
  Given any reflexive globular set $G$ there exists, for any judgement
  $\mathcal{J}$ of $\mathbb{T}_{\kappa}[G]$, a skeletal judgement
  $\textnormal{sk}(\mathcal{J})$ together with a substitution
  $\sigma_{\mathcal{J}}$ consisting entirely of basic terms such that 
  \begin{align*}
    \textnormal{sk}(\mathcal{J})[\sigma_{\mathcal{J}}] & \equiv\mathcal{J}.
  \end{align*}
  \begin{proof}
    Fix a derivation of $\mathcal{J}$.  There is a corresponding
    finite reflexive globular set $H$ obtained as the sub-reflexive
    globular set of $G$ generated by those cells of $G$ occurring in
    the fixed derivation of $\mathcal{J}$.  By the observations above
    we have $\textnormal{sk}_{H}(\mathcal{J})$ satisfying
    (\ref{eq:skeleton}).  But $\textnormal{sk}_{H}(\mathcal{J})$ is
    also a derivable term in $\mathbb{T}_{\kappa}[G]$ and $\sigma_{H}$
    is also a substitution in $\mathbb{T}_{\kappa}[G]$.  As such, we
    may take $\textnormal{sk}(\mathcal{J})$ to be
    $\textnormal{sk}_{H}(\mathcal{J})$ and $\sigma_{\mathcal{J}}$ to
    be $\sigma_{H}$.
  \end{proof}
\end{lemma}
Note that there may exist more than one skeletal judgement
$\textnormal{sk}(\mathcal{J})$ and more than one
$\sigma_{\mathcal{J}}$ satisfying the equation from Lemma
\ref{lemma:skeleta} since one may have multiple derivations of the
same judgement $\mathcal{J}$ which employ different basic terms.  We
might hope to choose a ``minimal'' derivation in some way and define a
canonical skeleton of $\mathcal{J}$ in that way, but this is also not
possible in the presence of truncation rules.  For example, if
$p,q:\idn{A}{n}(a,b)$ and $t:B(a)$ where we are in the situation that
trunaction applies at level $n$, then we have two distinct derivations
of $t:B(b)$ which employ distinct contexts of basic terms.
Nonetheless, any two derivations which employ the same basic terms
will give rise to the same skeletal judgement and substitution.

\begin{remark}
Instead of the above model theoretic proof that for every derivable judgement there exists 
a derivable skeletal judgement and a suitable substitution to recover the original judgement
one can also work purely syntactically by defining a translation on the level of 
expressions which replaces every basic term $\ulcorner a \urcorner$\/ by 
a variable $v_a$\/ and which commutes with all other formation rules. Then
one can show by a straightforward induction on derivations that if all basic terms in a derivation
of a judgement come from the finite globular set $H$, then the translated judgement is derivable
when we work in the context $\Delta_H$.
\end{remark}

One of the advantages of having Lemma \ref{lemma:skeleta} at our
disposal is that it allows us to give an alternative characterization of the monad
$M_{\kappa}$.  Let $\tilde{M}_{\kappa}(G)$ be the reflexive globular
set which has as cells (of the appropriate level) equivalence classes
of tuples $(\Gamma,\varphi,\sigma)$ such that
$\judge{\Gamma}{\varphi}$ is a skeletal judgement and $\sigma$ is a
substitution $()\to\Gamma$ such that $\varphi$ is required to have the
appropriate type (i.e., for $0$-cells
$\judge{\Gamma}{\varphi:\namefor{G}}$, et cetera).  Here 
\begin{align*}
  (\Gamma,\varphi,\sigma) \approx (\Delta,\psi,\tau) & \text{ if and
    only if }\judge{}{\varphi[\sigma]=\psi[\tau]}.
\end{align*}
Equivalently, by Lemma \ref{lemma:skeleta} (taking a skeleton of the
judgement $\judge{}{\varphi[\sigma]=\psi[\tau]}$),
$(\Gamma,\varphi,\sigma)\approx(\Delta,\psi,\tau)$ if and only if
there exists a skeletal context $\Theta$ extending both $\Gamma$ and
$\Delta$ and a substitution $\vartheta:()\to\Theta$ extending both $\sigma$ and
$\tau$ such that $\judge{\Theta}{\varphi=\psi}$.  $\tilde{M}_{\kappa}$
is readily seen to constitute a monad on the category of reflexive
globular sets using the same approach as in the definition of the
monad structure on $M_{\kappa}$.  Furthermore, there is an isomorphism
$\lambda:M_{\kappa}\iso\tilde{M}_{\kappa}$ of monads described as
follows.  Given $t$ in $M_{\kappa}(G)$ choose, by Lemma
\ref{lemma:skeleta}, a skeleton $(\Delta,\varphi,\sigma)$ of $t$ and
let $\lambda_{G}(t):=[\Delta,\varphi,\sigma]$.  This is 
independent of choice of representative from the definitional equality
class of $t$ and of the choice of skeleton by definition of
$\approx$.  Going the other way, given $[\Delta,\varphi,\sigma]$ let
$\lambda_{G}^{-1}$ send this data to the definitional equality class of
$\varphi[\sigma]$.  Again, this is trivially independent of the choice
of representative.  It is clear that $\lambda$ and $\lambda^{-1}$
constitute a natural isomorphism and that they are compatible with the
respective monad structures.  Henceforth we will freely employ this
isomorphism without explicit mention where convenient.


\subsection{Limits and colimits of algebras}\label{section:limits}

The aim of this section is to show that the monads $M_{\kappa}$ are
finitary.  One consequence of this is that $\mlcxs_{\kappa}$ is
cocomplete as well as being complete.

Assume given a filtered category $\II$ and a functor
$A:\II\to\rglob$.  Denote by $A^{\infty}$ the colimit of this
functor.  By definition, an $n$-cell of $A^{\infty}$ is an equivalence
class $[a]$ of $n$-cells of the coproduct $\coprod_{i}A(i)$, where
$a\in A(i)$ is equivalent to $a'\in A(j)$ if and only if there exist
arrows $\varphi:i\to k$ and $\varphi':j\to k$ in $\II$ such that
\begin{align*}
  A(\varphi)(a)& = A(\varphi')(a').
\end{align*}
We would like to prove that
\begin{align}\label{eq:filtered_colim}
  M_{\kappa}(A^{\infty}) & \iso \colimit_{i}M_{\kappa}\bigl(A(i)\bigr).
\end{align}
This will require an analysis of the valid derivations of
the theory $\mathbb{T}_{\kappa}[A^{\infty}]$.  To begin with, note that
$\mathbb{T}_{\kappa}[A^{\infty}]$ is obtained by augmenting
$\mathbb{T}_{\kappa}$ with the new basic type $\namefor{A^{\infty}}$ as well as with
basic terms $\namefor{[a]}$ of the appropriate types as described in Section
\ref{sec:monad}.

Assume we are given a derivable judgement $\mathcal{J}$ of the theory
$\mathbb{T}_{\kappa}[A^{\infty}]$.  Then it follows from Lemma
\ref{lemma:skeleta} that we have a skeleton
$\textnormal{sk}(\mathcal{J})$ and the corresponding substitution
$\sigma_{\mathcal{J}}$.  Suppose the basic terms occurring in
$\sigma_{\mathcal{J}}$ are $\namefor{[a_{1}]},\ldots,\namefor{[a_{n}]}$, where it is
possible that $n=0$, then it follows from the fact that $\mathcal{I}$
is filtered that we may find representatives $a_{1}',\ldots,a_{n}'$ of
the equivalence classes $[a_{1}],\ldots,[a_{n}]$ such that
$a_{1}',\ldots,a_{n}'$ are all in the same $A_{i}$ for some
$i\in\mathcal{I}$.  Therefore, we obtain a new substitution
$\sigma_{\mathcal{J}}'$ by substituting
$\namefor{a_{1}'},\ldots,\namefor{a_{n}'}$ instead of the
corresponding terms in $\sigma_{\mathcal{J}}$.  So the judgement
$\textnormal{sk}(\mathcal{J})[\sigma_{\mathcal{J}}']$ is derivable in
$\mathbb{T}_{\kappa}[A_{i}]$.  By considering the case where
$\mathcal{J}$ is an appropriate term judgement we are able to use this
line of reasoning to show that (\ref{eq:filtered_colim}) holds.
\begin{lemma}
  \label{lemma:filtered_colimits}
  Given a filtered category $\II$ together with a functor
  $A:\II\to\rglob$, there is an isomorphism (\ref{eq:filtered_colim})
  of reflexive globular sets.
  \begin{proof}
    Assume given a reflexive globular set $X$ together with a cocone
    $x_{i}:M_{\kappa}A(i)\to X$.  We now describe the induced map $\xi:M_{\kappa}A^{\infty}\to
    X$.  Given an element $\tau$ of $M_{\kappa}A^{\infty}$ we have, by
    the reasoning above, the term
    $\tau':=\textnormal{sk}(\tau)[\sigma_{\tau}']$ in
    $M_{\kappa}A^{i}$.  Therefore, we define:
    \begin{align*}
      \xi(\tau) & := x_{i}(\tau').
    \end{align*}
    This is immediately seen to be independent of the choice of
    $a_{1}',\ldots,a_{n}'$.  To see that the definition does not depend on the choice of
    skeleton we use the fact that two skeleta $\tau$ can be obtained
    as restrictions of another skeleton with a larger ambient
    context.  Moreover, it follows from this
    definition that each $x_{i}$ can be recovered
    by precomposing $\xi$ with the map $M_{\kappa}A(i)\to
    M_{\kappa}A^{\infty}$.  Finally, for uniqueness of $\xi$, observe that
    implicit the construction of $\xi$ above we have proved that
    for each cell $\tau$ of $M_{\kappa}A^{\infty}$ there exists some $i$
    such that $\tau$ is in the image of the map $M_{\kappa}A(i)\to M_{\kappa}A^{\infty}$.
  \end{proof}
\end{lemma}
By general results from category theory (see e.g.~\cite{BarrM:TTT}) we have
the following proposition:
\begin{proposition}
  \label{prop:complete_cocomplete}
  For each $\kappa=0,1,\ldots,\omega$, the category $\mlcxs_{\kappa}$
  is complete and cocomplete.
\end{proposition}


\section{Doppelg\"{a}ngers, $M_{0}$-algebras and $M_1$-algebras}\label{sec:doppelgangers}

Our purpose in this section is to characterize the category
$\mlcxs_0$\/ of algebras for the monad $M_0$\/ by proving that it is
equivalent to the category
of sets and to introduce the basic machinery which will allow us, in
Section \ref{sec:model} below, to characterize the category
$\mlcxs_{1}$. 
The 0-dimensional case is already instructive and provides
us with an opportunity to introduce some ideas and concepts which will
be put to work in a more complicated setting in the 1-dimensional case. 

We begin by discussing the reason why the results are nontrivial by
explaining the various ways in which the type theory
$\mathbb{T}_0[G]$\/ 
proves the existence of infinitely many duplicates of all of the
vertices, edges, and higher edges of the globular set $G$. These duplicates (here called
doppelg\"angers) must all be shown to be propositionally equal to elements of the original globular set $G$. 

Next, we establish the characterization of the $M_{0}$-algebras in a number of steps, making
use of the set-theoretic interpretation of extensional type theories
and the realizability semantics from \cite{Hofstra:CRMTT}. We will
concentrate on stating the main concepts and theorems and omit some of the detailed proofs, 
allowing the reader to follow the line of argument.

\subsection{Doppelg\"angers}\label{subsec:doppelgangers}

Fix a globular set $G$\/ and consider the type theory $\mathbb{T}_\omega[G]$\/ (or any extension of it). It is clear from the definition of the theory
$\mathbb{T}_\omega[G]$\/ that every vertex $a \in G_0$\/ is represented as a term in $\mathbb{T}_\omega[G]$, namely $a:G$. (We shall, as before, not
distinguish between an actual element in $G$\/ and its ``name" in the type theory.) Similarly, every
1-dimensional edge $f \in G_1$\/ is represented by $f:\idn{G}{}(a,b)$, where $s(f)=a, t(f)=b$, and so on in higher dimensions. One might, at first sight,
conjecture that these are the only judgements of this form, i.e. that whenever $\mathbb{T}_\omega[G]$\/ derives $\tau:G$\/ for a closed term $\tau$,
then $\tau$\/ must be an element of $G$\/ already, and whenever $\mathbb{T}_\omega[G]$\/ derives $\sigma:\idn{G}{}(a,b)$\/ then $\sigma \in G_1$\/ already.
However, things are more complicated than that, due to the elimination
rule for identity types. 

Suppose, for example, that we have $a,b,c \in G_0$\/ and a
non-reflexivity term \mbox{$f:a \to b$}\/ in $G_1$. Now we can consider the following derivation:
\begin{prooftree}
  \AxiomC{$\judge{x:G,y:G,z:\idn{G}{}(x,y)}{G:\type}$}
  \noLine
  \UnaryInfC{$\judge{x:G}{c:G}$}
  \noLine
  \UnaryInfC{$\judge{}{f:\idn{G}{}(a,b)}$}
  \RightLabel{$\id{}$ elimination}
  \UnaryInfC{$\judge{}{\Jelim\bigl([x:G]c,a,b,f):G}$}
\end{prooftree}
This creates a new term of type $G$\/ which we denote by $c\langle f\rangle$; we call it the \myemph{doppelg\"anger} of $c$\/ (at $f$). This term is not definitionally
equal to any of $a,b,c$. However, it is propositionally equal to $c$: this we can see from the derivation
\begin{prooftree}
  \AxiomC{$\judge{x:G,y:G,z:\idn{G}{}(x,y)}{\idn{G}{}(c,\Jelim\bigl([v:G]c,x,y,z)):\type}$}
  \noLine
  \UnaryInfC{$\judge{x:G}{\rr(c):\idn{G}{}(c,\Jelim\bigl([v:G]c,x,x,\rr(x))}$}
  \noLine
  \UnaryInfC{$\judge{}{f:\idn{G}{}(a,b)}$}
  \RightLabel{$\id{}$ elimination}
  \UnaryInfC{$\judge{}{\Jelim\bigl([x:G]c,a,b,f):\idn{G}{}(c,\Jelim\bigl([v:G]c,a,b,f))}$}
\end{prooftree}
showing that there is a term witnessing $c \simeq c\langle f\rangle$\/ (note that by the conversion rule the second premise reduces to 
$\judge{x:G}{\rr(c):\idn{G}{}(c,c)}$\/ so that the trivial term is well-defined).

Of course, this idea works in general: given any term $\tau:T$\/ and any (non-reflexivity) identity proof $f:\idn{A}{}(a,b)$\/ we may form 
\[ \tau\langle f\rangle := \Jelim([x:A]\tau,a,b,f):T\]
and then show that $\tau \simeq \tau\langle f\rangle$.

There are other ways to create doppelg\"angers: consider again $f:a \to b \in G_1$\/ and form 
\[f^\sharp:=\Jelim_{[x,y:G,z:\idn{G}{}(x,y)]G}([x:G]x,a,b,f):G.\] 
This term is a new vertex which is homotopic to both $a$\/ and $b$\/ (again this is proved by defining a suitable witness using the J-rule).

Yet another possibility is to construct 
\[ f^\flat:=\Jelim_{[x,y:G,z:\idn{G}{}(x,y)]\idn{G}{}(x,y)}([x:G]\rr(x),a,b,f):\idn{G}{}(a,b),\] 
which turns out to be homotopic to $f$.

While in the above examples of doppelg\"angers it is easy to show that each of the newly created terms is, up to homotopy, equal to a basic term coming from
the original globular set, it is not clear why this would always be the case, i.e. why for \emph{every} term derivable in $\mathbb{T}_\omega[G]$\/ there is
a suitable homotopy. Moreover, it will be seen in the next section that the elimination rule for identity types does in certain instances give genuinely
new terms which are not homotopic to any basic term (namely, the formal composites which are used to give the Martin-L\"of complexes their categorical
structure).

\subsection{$M_{0}$-algebras}\label{subsec:algebras_t0}

We now study the category of algebras $\mlcxs_0$\/ for the monad $M_0$. We fix a reflexive globular set $G$, and consider $M_0(G)$, the free
algebra on $G$. 

\begin{lemma}
The reflexive globular set $M_0(G)$\/ is constant. 
\end{lemma}

\begin{proof}
Since the theory $\mathbb{T}_0[G]$\/ satisfies the reflection rule, it follows that any term $\tau:\idn{G}{n}(a,b)$\/ is definitionally equal to 
a reflexivity term (see Subsection \ref{subsec:truncation} above). 
Hence for $n>0$, the elements of $M_0(G)_n$\/ are all degenerate, and the globular set $M_0(G)$\/ is completely
determined by its vertices.
\end{proof}

Thus in order to characterize the globular set $M_0(G)$, it suffices to understand the set $M_0(G)_0$\/ of its vertices. Recall from the construction 
of the monad $M_0$\/ that the elements of $M_0(G)_0$\/ are equivalence classes of closed terms $\tau:G$, where two of these are identified if the theory proves
that they are definitionally equal.
We begin by noting that there is a canonical map from $\pi_0(G)$\/ to $M_0(G)_0$, induced by the coequalizer
\begin{align}\label{eqn:pi_0}
\xymatrix{
G_1 \ar@<.5ex>[r]^s \ar@<-.5ex>[r]_t & G_0 \ar@{->>}[r]^e \ar[dr]_{\eta_0} & {\pi_0(G)} \ar@{-->}[d]^p \\
& & M_0(G)_0 }
\end{align}

Here, the map $\eta_0$\/ is the component of the unit $\eta:G \to M_0(G)$\/ at dimension 0. For every $f \in G_1$\/ with $s(f)=a, t(f)=b$\/ 
there is an axiom $f:\idn{G}{}(a,b)$\/ in $\mathbb{T}_0[G]$; by
truncation this forces $a=b$\/ in the theory, and hence $a$\/ and $b$\/ 
are identified as well in $M_0(G)$; hence $\eta_0s=\eta_0t$.

We would like to show that $p$\/ is a bijection; this would prove that $M_0$\/ is isomorphic to the (idempotent) monad $\Delta\pi_0$\/ on $\rgsets$, and in 
particular it would follow that the category of $M_0$-algebras is just the category of sets.

The first step in proving this is to exploit the fact that extensional ML type theories may be modelled in locally cartesian closed categories (see the original work
 of Seely~\cite{Seely:LCCCTT}, or the expository texts~\cite{Jacobs:CLTT, Hofmann:SSDT}). In particular, these theories may be soundly interpreted in the category of sets. 
 More concretely, the theory $\mathbb{T}_0[G]$\/ has the following set-theoretic model: interpret the basic type $G$\/ by the set $\pi_0(G)$, and interpret the
 basic terms $a:G$\/ by the element $[a]$, the connected component of $a$\/ in $G$.

\begin{lemma}\label{lemma:set_interpretation}
The above interpretation extends to a model of $\mathbb{T}_0[G]$\/ in the category of sets.
\end{lemma}

\begin{proof}
We need only verify the new axioms and the new conversion rule of the theory; if these are valid under the interpretation then the result follows by soundness.
By construction, the judgements $a:G$\/ for $a \in G$\/ are valid. The identity types $\idn{G}{}(a,b)$\/ will be interpreted in a degenerate way,
namely as the emptyset when $[a]\neq[b]$ and as the one element set
when $[a]=[b]$. Thus if we have an element $f \in G_1$\/ with $s(f)=a$\/ and $t(f)=b$, then the interpretation of $f$\/ may be taken to be $[a]=[b]$,
since the reflection rule allows us to derive $a=b$\/ from the axiom $f:\idn{G}{}(a,b)$. Similar reasoning works to show that the term judgements associated
to higher cells of $G$\/ are soundly interpreted. Finally, the new conversion rule $i(a)=\rr(a)$\/ holds under the interpretation since both sides of the equation
will be interpreted as $[a]$.
\end{proof}

The soundness of this interpretation guarantees that the map $p$\/ is injective: indeed, given two connected components $[a]$\/ and $[b]$\/ of $G$, suppose
that $p[a]=p[b]$. Then $\mathbb{T}_{0}[G]$\/ proves that $a=b$. But then this
equation should hold in the model $\pi_0(G)$, i.e. $[a]=[b]$\/ as
elements of $\pi_0(G)$.  In particular, the interpretation yields a
map $q:M_{0}G\to\pi_{0}(G)$ of reflexive globular sets which is a
retract of $p$.  We will show in the next section that $q$\/ is in fact the inverse of $p$.

\subsection{Combinatorial realizability models}\label{subsec:CRM}

In \cite{Hofstra:CRMTT} it is shown how to construct models of type 
theories such as $\mathbb{T}_{\kappa}[G]$ in such a way that the
interpretations of terms will provide additional data regarding the
syntax of these theories.  These models are called
\myemph{combinatorial realizability models} and can be seen as a
generalized form of realizability model in the usual sense where the
realizers of terms can be, intuitively, some kind of
combinatorial data (in the cases we care about they will usually be
edges constructed in the
syntax of the theory). We will explain the conditions required in
order for such a model to exist, but the proof of this fact is
somewhat involved and can be found in \cite{Hofstra:CRMTT}.
\begin{definition}
  Given a reflexive globular set $G$, a \myemph{notion of
    1-realizability for $G$} is a functor
  $\textnormal{real}:\Pi_{1}(M_{1}(G))\to\sets$.  
(Here, $\Pi_1(X,\alpha)$\/ is the underlying ``fundamental" groupoid of the complex 
$(X,\alpha)$, see next section 
for details.) Similarly, a notion
  of 0-realizability for $G$ is a functor $\pi_{0}(M_{0}(G))\to\sets$
\end{definition}
We will often write $\tau\forces t:\namefor{G}$ to indicate that $t$
is an element of $\textnormal{real}(t)$ and, given $f:t\to s$ in
$\Pi_{1}(M_{1}(G))$ we write $\tau\cdot f$ for
$\textnormal{real}(f)(\tau)$.  The theorem regarding combinatorial
realizability models from \cite{Hofstra:CRMTT} can then be stated
precisely as follows:
\begin{theorem}[\cite{Hofstra:CRMTT}]\label{theorem:CRM}
  Given a notion of 1-realizability for $G$ satisfying the following
  conditions:
  \begin{itemize}
  \item For each vertex $a$ of $G$, there exists a realizer
    $\alpha_{a}\forces \namefor{a}:\namefor{G}$.
  \item For each edge $f:a\to b$ in $G$, we have $\alpha_{a}\cdot
    f=\alpha_{b}$.
  \end{itemize}
  there exists a sound and complete model of $\mathbb{T}_{1}[G]$ in
  which a closed term $t$ of type $\namefor{G}$ is
  interpreted as a realizer $\alpha_{t}\forces t:\namefor{G}$ and in which a
  closed term $f:\idn{\namefor{G}}{}(t,s)$ is interpreted as a proof
  that $\alpha_{t}\cdot f=\alpha_{s}$.
\end{theorem}

The following application of these models will be used twice in the remainder of the paper:

\begin{theorem}\label{thm:nat} (\cite{Hofstra:CRMTT})
Let $G$\/ be a graph, $H$\/ a groupoid and let $P,Q:M_1G \to H$\/ be two functors. 
Suppose furthermore that we are given
a morphism $\alpha_a: P(a) \to Q(a)$\/ for each basic term $a$
such that, for each basic term $f:a \to b$, the following diagram commutes:
\[
\xymatrix{
P(a) \ar[d]_{P(f)} \ar[r]^{\alpha_a} & Q(a) \ar[d]^{Q(f)} \\
P(b) \ar[r]_{\alpha_b} & Q(b).
}
\]
Then there exists a natural transformation $\alpha:P \Rightarrow Q$\/ whose component 
at a basic term $a$\/ is $\alpha_a$.
\end{theorem} 

There is a 0-dimensional variant of the model construction which allows us to construct
sound and complete models for theories of the form $\mathbb{T}_0[G]$. In this case,
a notion of 0-realizability is a Set-valued functor on $M_0(G)$, regarded as a discrete groupoid.

\begin{corollary}\label{cor:0real}
  Given a notion of 0-realizability for $G$\/ satisfying the first condition of Theorem~\ref{theorem:CRM}, 
there exists a sound and complete model
  of $\mathbb{T}_{0}[G]$ in which a closed term $t$ of type
  $\namefor{G}$ is interpreted as a realizer $\alpha_{t}\forces t:\namefor{G}$.
\end{corollary}
It follows from Corollary \ref{cor:0real} that the map
$q:M_{0}G\to\pi_{0}G$ is the inverse of $p:\pi_{0}G\to M_{0}G$ (both
of these are described in Section \ref{subsec:algebras_t0} above).
\begin{lemma}
  \label{lemma:M0pi0}
  For any $G$, the map $p:\pi_{0}G\to M_{0}G$ is an isomorphism with
  inverse the map $q$ described in Section \ref{subsec:algebras_t0}.
  \begin{proof}
    Denote by $\bar{\cdot}:M_{0}G\to M_{0}G$ the function $p\circ q$.
    Let $\tau\forces t:\namefor{G}$ if and only if
    $\tau:\idn{\namefor{G}}{}(t,\bar{t})$.  This is a notion of
    $0$-realizability.  Moreover, the hypotheses of Corollary
    \ref{cor:0real} are satisfied since
    $\rr(\namefor{a})\forces\namefor{a}:\namefor{G}$.  Therefore, this
    determines a combinatorial realizability model.  Let us denote by
    $\alpha_{t}$ the interpretation of $t:\namefor{G}$.  Then we have,
    for any vertex $t$ in $M_{0}G$,
    $\alpha_{t}:\namefor{G}(t,\bar{t})$ and so, by 0-truncation,
    $t=\bar{t}$, as required.
  \end{proof}
\end{lemma}
\begin{proposition}
  There is an isomorphism of categories $\mlcxs_{0}\iso\sets$.
  \begin{proof}
    By Lemma \ref{lemma:M0pi0} it follows that $M_{0}$ is isomorphic
    to the monad $\Delta\pi_{0}$ and therefore the resulting
    categories of algebras are isomorphic.
  \end{proof}
\end{proposition}

\subsection{Fundamental groupoids of $M_{1}$-algebras}\label{subsec:T1_are_groupoids}

The aim of this section is to generalize the basic setup from Section
\ref{subsec:algebras_t0} to the case of $M_{1}$-algebras and the
category $\mlcxs_{1}$.  Contrary to what one might expect, this category 
is not equivalent to the category of groupoids. However, there is an
adjunction

\begin{align}\label{eq:ad_mlcx_1_grpd}
  \xy
  {\ar@/_.5cm/(0,0)*+{\mlcxs_{1}};(30,0)*+{\groupoids}};
  {\ar@/^.5cm/@{<-}(0,0)*+{\mlcxs_{1}};(30,0)*+{\groupoids}};
  {(15,0)*+{\perp}};
  \endxy
\end{align}
analogous to the adjunction between topological spaces or, better yet,
homotopy $1$-types, and groupoids.  

We first describe the right adjoint $\Pi_{1}:\mlcxs_{1}\to\groupoids$
which allows us to regard $M_{1}$-algebras as groupoids (where
groupoids are themselves as reflexive
globular sets in which all $n$-cells are degenerate for $n\geq 2$
provided with the additional structure of composites and inverses).
That a $M_{1}$-algebra can be endowed with the structure of a groupoid
follows immediately from the construction of composition
and inverse operations --- as well as the corresponding propositional
equalities witnessing the associativity, unit and inverse laws --- by Hofmann and Streicher
\cite{Hofmann:GITT}.  However, we will later require some of the details of
the proof of this fact and we therefore describe the construction explicitly.  First, recall that, given
any type $A$ together with terms $a,b:A$ and
$f:\idn{A}{}(a,b)$, the \myemph{inverse}
$f^{-1}:\idn{A}{}(b,a)$ of $f$ is defined to be the following elimination
term:
\begin{align*}
  f^{-1} & := \Jelim_{[x,y:A,z:\idn{A}{}(x,y)]\idn{A}{}(y,x)}\bigl([x:A]\rr(x),a,b,f\bigr).
\end{align*}
Moreover, when there exists a further propositional equality
$g:\idn{A}{}(b,c)$, the \myemph{composite} $(g\cdot f)$ \myemph{of $g$ with
  $f$} is defined to be the term
$\app\bigl(\Jelim(\lambda_{v}.v,b,c,g),f\bigr)$, where the $J$-term here is
written in full as
\begin{align*}
  \Jelim_{[x,y:A,z:\idn{A}{}(x,y)]\idn{A}{}(a,y)^{\idn{A}{}(a,x)}}\bigl([x:A]\lambda_{v:\idn{A}{}(a,x)}.v,b,c,g\bigr):\idn{A}{}(a,c)^{\idn{A}{}(a,b)}.
\end{align*}
We will use these operations on terms of identity type to define the
composition and inverses for $M_{1}$-algebras.  To this end, let an object $G$ of $\mlcxs_{1}$ be
given with action $\gamma:M_{1}(G)\to G$.  Of course, we will regard
$G$ as a groupoid with objects the vertices of $G$ and arrows the
edges of $G$.  Identities are given by the edges of the form $i(a)$
for $a$ a vertex.  In order to define composition in $G$ let a
composable pair of edges $f,g$ in $G$ be given with
\begin{align*}
  \xy
  {\ar^{f}(0,0)*+{a};(20,0)*+{b}};
  {\ar^{g}(20,0)*+{b};(40,0)*+{c.}};
  \endxy
\end{align*}
By definition, both of these edges (and their endpoints) are
represented by corresponding terms $f:\idn{G}{}(a,b)$ and $g:\idn{G}{}(b,c)$ in the
theory $\mathbb{T}_{1}[G]$.  As such, the composite $(g\cdot
f):\idn{G}{}(a,c)$, as defined above, exists and we define the result of composing
$f$ with $g$ in $G$ to be the edge obtained by appling the action of
$G$ to $(g\cdot f)$.  I.e.,
\begin{align*}
  (g\circ f) & := \gamma\bigl(\namefor{a},\namefor{c};\;(\namefor{g}\cdot \namefor{f})\bigr).
\end{align*}
This edge possesses the appropriate source and target since $\gamma$
is an arrow in $\rgsets$.  Likewise, the inverse $f^{-1}$ of $f$ is
defined by setting
\begin{align*}
  f^{-1} & := \gamma(\namefor{b},\namefor{a};\;\namefor{f}^{-1}),
\end{align*}
where $f^{-1}$ on the right-hand side is the inverse of the term
$\namefor{f}$, as defined above.  

With these definitions, the groupoid laws are a consequence of their up-to propositional
equality counterparts (for which see \cite{Hofmann:GITT}) together
with the $1$-truncation rule.  In this way the unit law is an
immediate consequence of the fact that $\rr(\namefor{a})=\namefor{i(a)}$.  For the associativity law,
suppose we are given $f$ and $g$ as above together with a further edge
$h:c\to d$ in $G$.  To prove the associative law $h\circ (g\circ f)  =
(h\circ g)\circ f)$ holds it suffices to show that
\begin{align}\label{eq:M_1_associative}
  \gamma\biggl(\namefor{h}\cdot\namefor{\gamma\bigl(\namefor{g}\cdot\namefor{f}\bigr)}\biggr) &= \gamma\biggl(\namefor{\gamma\bigl(\namefor{h}\cdot\namefor{g}\bigr)}\cdot\namefor{f}\biggr),
\end{align}
where we have omitted all but the final entries of lists of terms as
the missing entries are evident in this case.  To see that this is
indeed the case observe that the left-hand side of (\ref{eq:M_1_associative}) is equal to
\begin{align*}
  \gamma\bigl(\gamma_{*}(\namefor{\namefor{h}})\cdot\gamma_{*}(\namefor{(\namefor{g}\cdot\namefor{f})})\bigr)
  & = \gamma\circ
  M_{1}(\gamma)(\namefor{\namefor{h}}\cdot\namefor{(\namefor{g}\cdot\namefor{f})})\\
  & =
  \gamma\circ\mu_{G}(\namefor{\namefor{h}}\cdot\namefor{(\namefor{g}\cdot\namefor{f})})\\
  & =\gamma\bigl(\namefor{h}\cdot(\namefor{g}\cdot\namefor{f})\bigr)
\end{align*}
where the penultimate equality is by the multiplication law for
actions.
By the remarks above, $\namefor{h}\cdot(\namefor{g}\cdot\namefor{f})$
is definitionally equal to
$(\namefor{h}\cdot(\namefor{g}\cdot\namefor{f})$.  A
dual calculation reveals that the right-hand side of
(\ref{eq:M_1_associative}) is equal to
$\gamma((\namefor{h}\cdot\namefor{g})\cdot\namefor{f})$.  Lastly, that
$f^{-1}$ is the inverse of $f$ is straightforward using similar
reasoning.  That is, we have described a groupoid $\Pi_{1}(G,\gamma)$
constructed from a $M_{1}$-algebra $(G,\gamma)$.

In slightly more abstract terms this construction can be described as
follows.  For $G$ an arbitrary reflexive globular set, let
$\freegpd(G)$ denote the free groupoid (regarding the free groupoid
monad as a monad on reflexive globular sets) on $G$.  Recall that
$\freegpd(G)$ has the same vertices as $G$, and arrows
$a\to b$ in $\freegpd(G)$ are a zig-zag paths
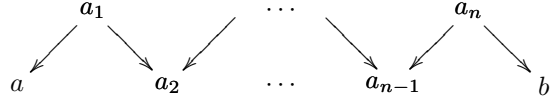
\begin{align*}
  \xy
  {\ar@{<-}(0,0)*+{a};(10,10)*+{a_{1}}};
  {\ar(10,10)*+{a_{1}};(20,0)*+{a_{2}}};
  {\ar@{<-}(20,0)*+{a_{2}};(30,10)*+{}};
  {\ar(40,10)*+{};(50,0)*+{a_{n-1}}};
  {\ar@{<-}(50,0)*+{a_{n-1}};(60,10)*+{a_{n}}};
  {\ar(60,10)*+{a_{n}};(70,0)*+{b}};
  {(35,10)*+{\cdots}};
  {(35,0)*+{\cdots}};
  \endxy
\end{align*}
of edges in $G$ modulo the evident relations forcing the groupoid laws
to hold.  There is then, for each $G$, a map $\Phi_{G}:\freegpd(G)\to
M_{1}(G)$ of globular sets which sends an equivalence class of such
``formal composites'' from $\freegpd(G)$
to the term representing the result of taking inverses and composites
of its edges using the type theoretic inverses and composites
described above.  These maps constitute a morphism of monads
$\freegpd\to M_{1}$ and therefore induce a functor
$\Pi_{1}:\mlcxs_{1}\to\groupoids$.  Explicitly, $\Pi_{1}(G,\gamma)$ is
given by the underlying globular set $G$ together with the action
$\gamma\circ\Phi_{G}:\freegpd(G)\to G$.  Moreover, $\Phi_{G}$ is
actually the canonical functor $\Phi_{G}:\freegpd(G)\to\Pi_{1}(M_{1}(G))$ extending 
the unit $G\to M_{1}(G)$ extends along the unit
$\eta'_{G}$ for $\freegpd$:
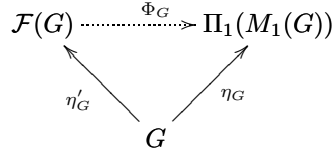
\begin{align*}
  \xy
  {\ar@{..>}^{\Phi_{G}}(0,15)*+{\freegpd(G)};(30,15)*+{\Pi_{1}(M_{1}(G))}};
  {\ar^{\eta'_{G}}(15,0)*+{G};(0,15)*+{\freegpd(G)}};
  {\ar_{\eta_{G}}(15,0)*+{G};(30,15)*+{\Pi_{1}(M_{1}(G))}};
  \endxy
\end{align*}
We sometimes call $\Pi_{1}(G,\gamma)$ the \emph{fundamental groupoid
  of $(G,\gamma)$}.  It follows from a general result of Kelly~\cite{Kelly:UTTC}
 (Theorem 25.4) that $\Pi_{1}$ possesses a left-adjoint
$K:\groupoids\to\mlcxs_{1}$.  We will return to a dicussion of this
adjunction later.  First we will turn to a proof that the maps
$\Phi_{G}$ constitute an equivalence of categories.

\subsection{Interpretation of $\mathbb{T}_{1}[G]$ using the free
  groupoid on $G$}\label{subsec:freegpd_interp}

The theory $\mathbb{T}_{1}[G]$ is soundly modelled using groupoids by extending
the interpretation from \cite{Hofmann:GITT} by the following
additional clauses:
\begin{itemize}
\item The new type $\namefor{G}$ is interpreted as the free groupoid on $G$:
  \begin{align*}
    \lscott\namefor{G}\rscott & := \freegpd(G).
  \end{align*}
\item The new terms basic $\namefor{a}$ of type $\namefor{G}$ are
  interpreted by the objects of $\freegpd(G)$ which they represent:
  \begin{align*}
    \lscott\namefor{a}\rscott & := a.
  \end{align*}
\item The new basic terms $\namefor{f}$ of identity type
  $\namefor{G}(\namefor{a},\namefor{b})$ are likewise interpreted as
  the arrows they represent
  \begin{align*}
    \lscott\namefor{f}\rscott & := f.
  \end{align*}
\item If $\namefor{\alpha}$ is a new basic term of type
  $\namefor{G}^{n}(\namefor{\alpha^{0}_{0}},\ldots,\namefor{\alpha^{n-1}_{1}})$, for $n>1$, then
  \begin{align*}
    \lscott\namefor{\alpha}\rscott & := \alpha^{1}_{0}.
  \end{align*}
\end{itemize}
With these definitions, the axioms of $\mathbb{T}_{1}[G]$ are clearly
satisfied.  We now remind the reader how the particular kinds of terms
we are interested in are interpreted in this model.  To begin with
recall that the identity type $\judge{x,y:\namefor{G}}{\namefor{G}(x,y):\type}$ is interpreted as the functor
$I_{G}:\freegpd(G)\times\freegpd(G)\to\groupoids$ which sends a pair of objects $(a,b)$ of $\freegpd(G)$
to the discrete groupoid $\freegpd(G)(a,b)$ and which sends an
arrow $(\alpha,\beta):(a,b)\to(a',b')$ to the functor
$\freegpd(G)(a,b)\to\freegpd(G)(a',b')$ with action $f\mapsto(\beta\circ
f\circ\alpha^{-1})$.  The extended context
$(x,y:\namefor{G},z:\namefor{G}(x,y))$ is interpreted as the result of
applying the Grothendieck construction $\int I_{G}$ to $I_{G}$.  In
this instance, $\int I_{G}$ coincides with the arrow category
$\freegpd(G)^{\to}$.  As such, the elimination data
$\judge{x:\namefor{G}}{\varphi(x):B(x,x,\rr(x))}$ is interpreted by
a functor $\lscott B\rscott:\freegpd(G)^{\to}\to\groupoids$ together
with a functor $\lscott\varphi\rscott:\freegpd(G)\to\int\lscott B\rscott$ such that
\begin{align*}
  \xy
  {\ar^{\varphi}(0,15)*+{\freegpd(G)};(25,15)*+{\int\lscott B\rscott}};
  {\ar^{\pi}(25,15)*+{\int\lscott B\rscott};(25,0)*+{\freegpd(G)^{\to}}};
  {\ar@/_.5pc/_{r}(0,15)*+{\freegpd(G)};(25,0)*+{\freegpd(G)^{\to}}};
  \endxy
\end{align*}
commutes.  I.e., for an object $a$ of $\freegpd(G)$, $\varphi(a)$ is a tuple
composed of $1_{a}:a\to a$ together with an object, which we denote by
$a_{\varphi}$, of the groupoid $\lscott B\rscott(1_{a}:a\to a)$.  For
an arrow $\alpha:a\to a'$ of $\freegpd(G)$, $\varphi(\alpha)$ is then
a tuple composed of $\alpha$ itself together with an arrow
\begin{align*}
  \lscott B\rscott\left(\vcenter{\xy
        {\ar^{1_{a}}(0,10)*+{a};(10,10)*+{a}};
        {\ar_{1_{a'}}(0,0)*+{a'};(10,0)*+{a'}};
        {\ar_{\alpha}(0,10)*+{a};(0,0)*+{a'}};
        {\ar^{\alpha}(10,10)*+{a};(10,0)*+{a'}};
        \endxy}\right)\vcenter{\xy
      {\ar^{\alpha_{\varphi}}(0,0)*+{(a_{\varphi})};(20,0)*+{a'_{\varphi}}};
      \endxy}
  \end{align*}
in the groupoid $\lscott B\rscott(1_{a'}:a'\to a')$.

The resulting elimination term
$\judge{x,y:\namefor{G},z:\namefor{G}(x,y)}{\Jelim(\varphi,x,y,z):B(x,y,z)}$
is interpreted as the section $J$ of the projection $\int\lscott
B\rscott\to\freegpd(G)^{\to}$ which sends an object $f:a\to b$ of $\freegpd(G)^{\to}$ to the pair consisting of $f$ and the object
\begin{align*}
  \lscott B\rscott\left(\vcenter{\xy
  {\ar^{1_{a}}(0,10)*+{a};(10,10)*+{a}};
  {\ar_{f}(0,0)*+{a};(10,0)*+{b}};
  {\ar_{1_{a}}(0,10)*+{a};(0,0)*+{a}};
  {\ar^{f}(10,10)*+{a};(10,0)*+{b}};
  \endxy}\right)\bigl(a_{\varphi}\bigr)
\end{align*}
of $\lscott B\rscott(a,b,f)$.  Similarly, the action of $J$ on an arrow 
\begin{align}\label{eq:ar_ar_cat}
  \vcenter{\xy
  {\ar^{f}(0,10)*+{a};(10,10)*+{b}};
  {\ar_{g}(0,0)*+{a'};(10,0)*+{b'}};
  {\ar_{\alpha}(0,10)*+{a};(0,0)*+{a'}};
  {\ar^{\beta}(10,10)*+{b};(10,0)*+{b'}};
  \endxy}
\end{align}
from $f:a\to b$ to $g:a'\to b'$ in $\freegpd(G)^{\to}$ is the pair
consisting of the arrow itself together with 
\begin{align*}
  \lscott B\rscott\left(\vcenter{\xy
      {\ar^{1_{a'}}(0,10)*+{a'};(10,10)*+{a'}};
      {\ar_{g}(0,0)*+{a'};(10,0)*+{b'}};
      {\ar_{1_{a'}}(0,10)*+{a'};(0,0)*+{a'}};
      {\ar^{g}(10,10)*+{a'};(10,0)*+{b'}};
      \endxy
    }\right)\bigl(\alpha_{\varphi}\bigr):
  \lscott B\rscott\left(\vcenter{\xy
      {\ar^{1_{a}}(0,10)*+{a};(10,10)*+{a}};
      {\ar_{g}(0,0)*+{a'};(10,0)*+{b'}};
      {\ar_{\alpha}(0,10)*+{a};(0,0)*+{a'}};
      {\ar^{\beta\circ f}(10,10)*+{a};(10,0)*+{b'}};
      \endxy
    }\right)
  \vcenter{\xy
    {\ar(0,0)*+{a_{\varphi}};(10,0)*+{}};
    \endxy
  }
    \lscott B\rscott\left(\vcenter{\xy
        {\ar^{1_{a'}}(0,10)*+{a'};(10,10)*+{a'}};
        {\ar_{g}(0,0)*+{a'};(10,0)*+{b'}};
        {\ar_{1_{a}}(0,10)*+{a'};(0,0)*+{a'}};
        {\ar^{g}(10,10)*+{a'};(10,0)*+{b'}};
        \endxy
      }\right)\bigl(a'_{\varphi}\bigr)
\end{align*}
So, for example, given a term $h:\namefor{G}(g, g')$
in $\mathbb{T}_{1}[G]$, consider $h^{-1}$.  The pattern type
$B(x,y,z)$ in this instance is $\idn{G}{}(y,x)$ and $\lscott B\rscott$ is the
functor sending $f:a\to b$ to the discrete groupoid $\idn{G}{}(b,a)$ and which
sends an arrow (\ref{eq:ar_ar_cat}) in $\freegpd(G)^{\to}$ to the
functor $\lambda_{v}.\alpha\circ v\circ\beta^{-1}:\idn{G}{}(b,a)\to
\idn{G}{}(b',a')$.  As such, it is straightforward to verify with the
description of the interpretation given above that
$\lscott h^{-1}\rscott$ is equal to the inverse $\lscott h\rscott^{-1}$ in
$\freegpd(G)$.  Similarly, given $f:\namefor{G}(a,b)$ and $g:\namefor{G}(b,c)$ in
$\mathbb{T}_{1}[G]$, it is straightforward to verify that the
interpretation commutes with composition in the sense that $\lscott
\namefor{g}\cdot\namefor{f}\rscott$ is equal
to $\lscott g\rscott\circ\lscott f\rscott$ in $\freegpd(G)$.  These
observations yield the following:
\begin{lemma}
  \label{lem:theta_prime}
  The assignment $\Psi_{G}:\Pi_{1}(M_{1}(G))\to\freegpd(G)$ which sends an $n$-cell
  $\vec{\alpha}=(\alpha^{0}_{0},\ldots,\alpha^{n-1}_{1},\alpha)$ of
  $M_{1}(G)$ to $\lscott\alpha\rscott$ is functorial.
  \begin{proof}
    By the results of Section \ref{subsec:T1_are_groupoids} it follows
    that $M_{1}(G)$ is a groupoid in with the result of composing
    1-cells $(a,b;f)$ and $(b,c;g)$ is $(a,c;g\cdot f)$.  Thus,
    because the interpretation function commutes with composition it
    follows that $\Psi_{G}$ is functorial (that $\Psi_{G}$ preserves
    identities is also straightforward).
  \end{proof}
\end{lemma}

\begin{theorem}
  \label{theorem:T1G_is_FG}
  Given a reflexive globular set $G$, $\Phi_{G}:\freegpd(G)\to
  \Pi_{1}(M_{1}(G))$ is an equivalence of categories.
  \begin{proof}
    It is an immediate consequence of the universal property of
    $\freegpd(G)$ that $\Phi_{G}$ is a section of $\Psi_{G}$.  
    Construct a combinatorial realizability model
    (as described in Section \ref{subsec:CRM} above) of
    $\mathbb{T}_{1}[G]$ where realizers of terms of type $\namefor{G}$
    are given by letting $\varphi\forces t:\namefor{G}$ if
    and only if $\varphi$ is a closed term
    \begin{align*}
      \judge{}{\varphi:\idn{\namefor{G}}{}\bigl(t,\Phi_{G}(\Psi_{G}(t))\bigr)}
    \end{align*}
    which satisfies $\Psi_{G}(\varphi)=1_{\psi_{G}(t)}$.  Realizers
    for basic terms $\namefor{a}$ are given by $\rr(\namefor{a})$.  It
    follows from the results of \cite{Hofstra:CRMTT} that this
    determines a model of type theory.  Let us denote the
    interpretation (i.e., realizer) of $t:\namefor{G}$ by
    $\alpha_{t}$.  Then, by virtue of the interpretations of identity
    types in the realizability model, it follows that these
    $\alpha_{t}$ give a natural isomorphism
    $\alpha:1_{\Pi_{1}(M_{1}(G))}\iso\Phi_{G}\circ\Psi_{G}$.
  \end{proof}
\end{theorem}
Theorem~\ref{theorem:main} shows that free $M_1$-algebras are, up to
equivalence, free groupoids. This might lead one to conjecture that
the category of $M_1$-algebras is equivalent to the category of
groupoids. However, that is not the case. 
The following example makes clear that different algebras may have the
same fundamental groupoid.
\begin{example}
Consider the following groupoid $G$: it has two objects $a$\/ and $b$,
exactly one arrow $f:a \to b$\/ and its inverse $g:b \to a$. We may
define a $M_1$-algebra
structure $\gamma:M_1(G) \to G$\/ on $G$\/ as follows: on objects,
$\gamma$\/ is defined by
\[ 
\gamma(v)=\left\{\begin{array}{ll}
a & \text{ if } \judge{}{v=\namefor{a}:G} \text{ is derivable } \\
b & \text{ otherwise}
\end{array} \right.
\]
Thus, all doppelg\"angers of vertices are sent to $b$. On 1-cells we define:
\[ \gamma(w)=
\left\{\begin{array}{ll}
f & \text{ if } \gamma(s(w))=a, \gamma(t(w))=b \\
g & \text{ if } \gamma(s(w))=b, \gamma(t(w))=a \\
1_a & \text{ if } \gamma(s(w))=a=\gamma(t(w)) \\
1_b & \text{ if } \gamma(s(w))=b=\gamma(t(w)) \\
\end{array} \right.\]

It is readily seen that this is a map of globular sets. To see that it
is a $M_1$-algebra, we remark that the unit law is trivially satisfied because
the algebra map sends any generator $\namefor{v}$\/ of $M_1(G)$\/ to
$v$. For the associativity law, consider an element $\tau$\/ of $M_1^2(G)_0$;
this is a term of the theory $\mathbb{T}_1[M_1(G)]$, which is
generated by basic terms of the form $\namefor{\sigma}$, where $\sigma$\/ is a term
of the theory $\mathbb{T}_1[G]$. Note that on the one hand
\[ (\gamma \circ M_1\gamma)(\tau)=a \Leftrightarrow M_1\gamma(\tau)=\namefor{a} \Leftrightarrow \tau=\namefor{\namefor{a}}, \]
while on the other hand
\[ (\gamma \circ \mu)(\tau)=a \Leftrightarrow \mu(\tau) = \namefor{a} \Leftrightarrow \tau=\namefor{\namefor{a}}, \]
showing that both maps agree in dimension 0. To show that they agree
in dimension 1 as well, one reasons in a similar fashion.

But clearly by symmetry there is another algebra structure on $G$,
call it $\delta$, defined by sending all doppelg\"angers to $a$\/ instead of $b$.
The identity map $G \to G$\/ is, however, not a map of
$M_1$-algebras. Indeed, any map of $M_1$-algebras commutes with the formation of
doppelg\"angers; for example, if $k$\/ is a map of algebras then
$M_1(k)$\/ must send the doppelg\"anger $a\langle f\rangle$\/ to
$k(a)\langle k(f)\rangle$, and
hence we must have $k\gamma(a\langle f \rangle) = \delta k(a)\langle
k(f)\rangle$, which is impossible if $k$\/ is the identity. 
For the same reason the only other possible map of groupoids, which
interchanges $a$\/ and $b$, cannot be a map of $M_1$-algebras.
\end{example}

Thus $M_1$-algebras carry more information than their fundamental
groupoids, and this information tells us how the formal composites and
doppelg\"angers are interpreted.
The fact that non-isomorphic algebras may have the same fundamental
groupoid is of course the analogue of the fact that non-homeomorphic
topological spaces may have the same fundamental groupoid.

In addition, the above example shows that $\Pi_{1}$\/ is not a full
functor. (However, it is easily seen to be faithful. )  Nonetheless, it will be
shown in Section \ref{sec:model} below that $K\adjoint\Pi_{1}$
constitutes a Quillen equivalence and it is to this that we now turn.


\section{The Quillen model structure on $\mlcxs_{1}$}\label{sec:model}

In this section we will only consider 1-truncated complexes, and we reduce clutter 
in the notation by dropping subscripts indicating this one-dimensionality.  Given an
object $(A,\alpha)$ of $\mlcxs_{1}$, we will sometimes denote the
composition in the resulting groupoid $\Pi_{1}(A,\alpha)$ by $\circ_{\alpha}$.

We begin by defining the three classes of morphisms for the model structure on $\mlcxs_1$:
\begin{description}
\item[Fibrations] a map $f$\/ of complexes is a \emph{fibration} when $\Pi_{1}(f)$ 
is an isofibration of groupoids.  We denote the class of fibrations by $\frak{F}$. 
\item[Weak Equivalences]a map $f$\/ of complexes is a \emph{weak equivalence} when $\Pi_{1}(f)$ 
is a weak categorical equivalence.  We denote by the class of weak equivalences by $\frak{W}$.
\item[Cofibrations] A map of complexes is a \emph{cofibration} when it has
 the left lifting property with respect to maps which
are simultaneously fibrations and weak equivalences. The class of cofibrations will be denoted by 
$\frak{C}$.
\end{description}



\subsection{Cotensor of complexes with graphs}

Let a  complex $(A,\alpha)$ be given together with a graph $X$.  We define a new
complex $(A^{X},\chi)$ as follows.  The underlying graph has as
0-cells graph homomorphisms $F:X\to A$ and as 1-cells natural
transformations.  Here naturality of a transformation $\eta:F\to G$
means that for each vertex $x$ of $X$ we have a 1-cell $\eta_{x}:Fx\to
Gx$ in $A$ such that, for $h:x\to y$ in $X$, we have
\begin{align*}
  \eta_{y}\circ_{\alpha}Fh & = Gh\circ_{\alpha}\eta_{x}.
\end{align*}
Now, fix a vertex $x$ in $X$.  We define an evaluation map
$\varepsilon_{x}:M_{1}(A^{X})\to M_{1}A$ as the map of expressions (trivially
seen to preserve derivable judgements) which sends $\namefor{F}$ to
$\namefor{Fx}$ and $\namefor{\alpha}$ to $\namefor{\alpha_{x}}$.
Before we can go any further we must make some observations regarding
these evaluation maps.  We begin with the following fact which follows
immediately from the definition of $\varepsilon_{x}$:
\begin{lemma}
  The evaluation map $\varepsilon_{x}:M_{1}(A^{X})\to M_{1}A$ is functorial.
\end{lemma}
For the following theorem we must construct a combinatorial
realizability model (see Section \ref{subsec:CRM}) of $\mathbb{T}_{1}[A^{X}]$:
\begin{theorem}
  Given an edge $f:x\to y$ in $X$, there is an induced natural
  transformation $\varepsilon_{f}:\varepsilon_{x}\to\varepsilon_{y}$.
  \begin{proof}
    We define $\varphi\forces t:A^{X}$ to hold if and only if
    $\varphi:\idn{A}{}(\varepsilon_{x}t,\varepsilon_{y}t)$.  This has
    a functorial action since given $g:\idn{A^{X}}{}(t,s)$ and
    $\varphi\forces t: A$ we have 
    \begin{align*}
      \varepsilon_{y}(g)\cdot\varphi\cdot\varepsilon_{x}(g)^{-1}:\idn{A}{}(\varepsilon_{x}s,\varepsilon_{y}s).
    \end{align*}
    Moreover, since $\varepsilon_{x}$ and $\varepsilon_{y}$ are
    functorial it follows that this action on realizers is also
    functorial.  Basic terms $\namefor{F}:A^{X}$ are realized by
    \begin{align*}
      \namefor{Ff}:\idn{B}{}(\namefor{Fx},\namefor{Fy}).
    \end{align*}
    With these definitions the conditions for a combinatorial
    realizability model are met and the existence of the natural
    transformation $\varepsilon_{f}$ follows (cf. Theorem~\eqref{thm:nat}).
  \end{proof}
\end{theorem}
We will write the component $\varepsilon_{x}(t)\to\varepsilon_{y}(t)$
of $\varepsilon_{f}$ at a vertex $t$ as $\varepsilon_{f}(t)$ and we
will assume that $\varepsilon_{1_{x}}$ is the identity.
We define the map $\chi:M_{1}(A^{X})\to A^{X}$ by
\begin{align*}
  \chi(t)(x) & := \alpha\bigl(\varepsilon_{x}(t)\bigr)
\end{align*}
for $t$ a vertex of $M_{1}(A^{X})$ and $x$ a vertex of $X$, and, for
$f:x\to y$ in $X$, we have 
\begin{align*}
  \chi(t)(f) & := \alpha\bigl(\varepsilon_{f}(t)\bigr).
\end{align*}
Next, for $g:t\to s$ in $M_{1}(A^{X})$ we define $\chi(t)\to\chi(s)$ by
taking at a vertex $x$ the map 
\begin{align*}
  \chi(g)_{x} & := \alpha\bigl(\varepsilon_{x}(g)\bigr).
\end{align*}
That this is a natural transformation is by naturality of the
$\varepsilon_{f}:\varepsilon_{x}\to\varepsilon_{y}$.  I.e., we have
proved the following:
\begin{lemma}
  The map $\chi$ is a graph homomorphism $M_{1}(A^{X})\to A^{X}$.
\end{lemma}
It now remains to show that this map gives $A^{X}$ a $M_{1}$-algebra structure.
\begin{lemma}
  $(A^{X},\chi)$ is a complex.
  \begin{proof}
    The unit law is trivial.  For the multiplication law assume given
    a term of $\mathbb{T}_{1}[M_{1}(A^{X})]$ of the form $\varphi(\namefor{\xi})$ where $\xi\in
    M_{1}(A^{X})$ and where $\varphi(-)$ is skeletal (note that we
    should really take an arbitrary list of basic terms
    $\namefor{\xi}$, but in the more general case the argument is
    identical to the one given here).  Then we must show that
    \begin{align*}
      \chi\bigl(\varphi(\xi)\bigr) & = \chi\bigl(\varphi(\namefor{\chi(\xi)})\bigr).
    \end{align*}
    It suffices to evaluate on $x\in X$.  We then have
    \begin{align*}
      \chi\bigl(\varphi(\namefor{\chi(\xi)})\bigr)(x) & =
      \alpha\bigl(\varepsilon_{x}(\varphi(\namefor{\chi(\xi)}))\bigr)\\
      & =
      \alpha\bigl(\varphi(\namefor{\alpha(\varepsilon_{x}(\xi))})\bigr)\\
      & = \alpha\bigl(\varphi(\varepsilon_{x}(\xi))\bigr)\\
      & = \alpha\bigl(\varepsilon_{x}(\varphi(\xi))\bigr),
    \end{align*}
    where the third equation is by the fact that $\alpha$ is an algebra
    and the fourth equation is by the fact that $\varphi(-)$ is skeletal.
  \end{proof}
\end{lemma}
We will denote by $X\pitchfork(A,\alpha)$ the algebra $(A^{X},\chi)$
when we do not want to have to mention the action $\chi$ and we call
this the \myemph{cotensor of $(A,\alpha)$ with $X$}.

Note that by construction of the map $\chi$, the diagram
\[
\xymatrix{
M(A^X) \ar[rr]^{M(ev_x)=\varepsilon_x} \ar[d]_\chi && MA \ar[d]^\alpha\\
A^X \ar[rr]_{ev_x} && A
}
\]
commutes, where $ev_x(F)=F(x)$. 
Thus the evaluation maps are actually algebra morphisms.

\subsection{The path object argument}

Where $\mathcal{F}$ denotes the free groupoid functor, we have the
following extremely useful fact:
\begin{lemma}
  \label{lemma:cotensor_Pi_1}
  For any graph $X$ and any complex $(A,\alpha)$, there is an
  isomorphism of groupoids
  \begin{align*}
    \Pi_{1}\bigl(X\pitchfork(A,\alpha)\bigr) & \iso \Pi_{1}(A,\alpha)^{\mathcal{F}X}.
  \end{align*}
  \begin{proof}
    This is routine using the universal property of
    $\mathcal{F}X$ and the definition of the edges in
    $X\pitchfork(A,\alpha)$ as natural transformations.
  \end{proof}
\end{lemma}
\begin{lemma}\label{lemma:path_objects}
  Each object $(A,\alpha)$ of $\mlcxs_{1}$ has a path object
  factorization.
  \begin{proof}
    Let $I$ be the graph with two vertices $0$ and $1$ and one
    non-trivial edge $0\to 1$.  Then $\mathcal{F}I$ is the usual
    ``interval'' $\mathbf{I}$ in the category of groupoids.  For any $(A,\alpha)$
    we have
    \begin{align}\label{eq:a_diagram}
      \xy
      {\ar^{r}(0,15)*+{(A,\alpha)};(30,15)*+{I\pitchfork(A,\alpha)}};
      {\ar_{\Delta}(0,15)*+{(A,\alpha)};(15,0)*+{(A,\alpha)\times(A,\alpha)}};
      {\ar^{p}(30,15)*+{I\pitchfork(A,\alpha)};(15,0)*+{(A,\alpha)\times(A,\alpha)}};
      \endxy
    \end{align}
    where 
    \begin{align*}
      r(a)(x) & := a
    \end{align*}
    and 
    \begin{align*}
      p(H) & := \langle H(0), H(1)\rangle.
    \end{align*}
    It is routine to verify that these are algebra homomorphisms.

    Using Lemma \ref{lemma:cotensor_Pi_1} and the fact that $\Pi_{1}$
    is a right-adjoint it follows that the result of applying
    $\Pi_{1}$ to (\ref{eq:a_diagram}) is 
    \begin{align*}
      \xy
      {\ar(0,15)*+{\Pi_{1}(A)};(30,15)*+{\Pi_{1}(A)^{\mathbf{I}}}};
      {\ar_{\Delta}(0,15)*+{\Pi_{1}(A)};(15,0)*+{\Pi_{1}(A)\times\Pi_{1}(A)}};
      {\ar(30,15)*+{\Pi_{1}(A)^{\mathbf{I}}};(15,0)*+{\Pi_{1}(A)\times\Pi_{1}(A)}};
      \endxy
    \end{align*}
    which, as is well known, constitutes a path object for
    $\Pi_{1}(A)$ in $\groupoids$.  Therefore, the original diagram
    (\ref{eq:a_diagram}) is a path object in $\mlcxs_{1}$.
  \end{proof}
\end{lemma}
Recall that Quillen's path object argument provides conditions under
which it is possible to transfer a model structure from a category
$\CC$ to a category $\DD$ along an adjunction $F\adjoint G$ for
$F:\CC\to\DD$.  
\begin{theorem}[Quillen]\label{theorem:Quillen}
  Assume given a cofibrantly generated model category $\CC$ together
  with an adjunction $F\adjoint G$ for $F:\CC\to\DD$ where $\DD$ is a
  cocomplete category with finite limits.  Assume furthermore that the
  following conditions are satisfied:
  \begin{enumerate}
  \item The left-adjoint $F$ preserves small objects.
  \item $\DD$ has a fibrant replacement functor.
  \item $\DD$ has a functorial path objects for fibrant objects.
  \end{enumerate}
  Then there is a model structure on $\DD$ in which a map is a
  fibration (weak equivalence) if and only if its image under $G$ is a
  fibration (weak equivalence) in $\CC$.
\end{theorem}
From this we are able to obtain the following:
\begin{theorem}\label{theorem:main}
  The definition of fibration and weak equivalence in $\mlcxs_{1}$
  given above determines a model structure on $\mlcxs_{1}$.
  \begin{proof}
    By Theorem \ref{theorem:Quillen}, Lemma \ref{lemma:path_objects}
    and the fact that every object in $\groupoids$ is fibrant (and hence that 
    every object of $\mlcxs_1$\/ is fibrant by definition), it
    suffices to prove that the left-adjoint
    $K:\groupoids\to\mlcxs_{1}$ of $\Pi_{1}$ preserves small objects.
    Note that both forgetful functors $\mlcxs_{1}\to\rglob$ and
    $\groupoids\to\rglob$ preserve and reflect filtered colimits, since
    both are finitarily monadic. 
    Therefore, the functor $\Pi_{1}:\mlcxs_{1}\to\groupoids$ must
    preserve filtered colimits, and in particular colimits of
    chains.  This last statement is equivalent to the preservation of
    small objects by $K$.
  \end{proof}
\end{theorem}

\subsection{The construction of the left adjoint $K$ of $\Pi_{1}$}

In order to prove that the adjunction $K\adjoint\Pi_{1}$ is a Quillen
equivalence it will be necessary to consider the transfinite construction of $K$
from~\cite{Kelly:UTTC} in more detail\footnote{The construction can also be seen as 
combining the transfinite construction for coequalizers in categories of algebras as detailed, for example,
in~\cite{BarrM:TTT} with the fact that the left adjoint to the functor $\Pi_1$\/ can be rendered as
a coequalizer.}.  Henceforth $G$ denotes a fixed groupoid.

Let us set $G(-1):=G$ and $G'(-1):=G$.  We construct $G(0)$ as
the following coequalizer (taken in the category of graphs):
\begin{align*}
  \xy
  {\ar^{M_{1}\gamma}(0,0)*+{M_{1}\mathcal{F}G};(30,0)*+{M_{1}G}};
  {\ar_{M_{1}\Phi_{G}}@/_1em/(0,0)*+{M_{1}\mathcal{F}G};(15,-10)*+{M_{1}^{2}G}};
  {\ar_{\mu_{G}}@/_1em/(15,-10)*+{M_{1}^{2}G};(30,0)*+{M_{1}G}};
  {\ar^{e(-1)}(30,0)*+{M_{1}G};(50,0)*+{G(0)}};
  \endxy
\end{align*}
where $\gamma$ is the action of $G$ (\emph{qua} groupoid) and
$\Phi_{G}$ is the canonical map induced by the groupoid structure on
$M_{1}G$.  Let $i(-1):G(-1)\to G(0)$ denote the composite
$e(-1)\circ\eta_{G(-1)}$.

In the next stages of the construction of $K$ we obtain $G(n+2)$ from
$G(n)$ and $G(n+1)$ as the following coequalizer
\begin{align*}
  \xy
  {\ar@/^1em/^{\mu_{G(n)}}(0,0)*+{M_{1}^{2}G(n)};(15,10)*+{M_{1}G(n)}};
  {\ar@/^1em/^{M_{1}\eta_{G(n)}}(15,10)*+{M_{1}G(n)};(30,0)*+{M_{1}^{2}G(n)}};
  {\ar_{1_{M_{1}^{2}G(n)}}(0,0)*+{M_{1}^{2}G(n)};(30,0)*+{M_{1}^{2}G(n)}};
  {\ar^{M_{1}e(n)}(30,0)*+{M_{1}^{2}G(n)};(60,0)*+{M_{1}G(n+1)}};
  {\ar^{e(n+1)}(60,0)*+{M_{1}G(n+1)};(90,0)*+{G(n+2)}};
  \endxy
\end{align*}
and we define $i(n+1):G(n+1)\to G(n+2)$ to be the composite
$e(n+1)\circ\eta_{G(n+1)}$.

$KG$ is then defined as the colimit (taken in the category of
reflexive globular sets) $\colimit_{n}G(n)$ of the diagram
consisting of the maps $i(n)$.  The action $\nu:M_{1}KG\to KG$ is the
canonical map induced by the maps $M_{1}G(n)\to G(n+1)\to KG$ together
with the fact that $M_{1}$ is finitary and hence preserves the colimit of the chain $G(i)$.

\subsection{The Quillen equivalence}

We now begin working towards the proof that the adjunction $K \dashv \Pi_1$\/ is a Quillen
equivalence. This will be done in several steps: first, with the aid of the groupoid semantics
we construct a sequence of functors $\lscott - \rscott_n:M_1G(n) \to G$\/ with suitable properties. 
Next, we construct, by induction on $n$, a realizability model of $\mathbb{T}_1[G(n)]$\/ which, 
using Theorem~\ref{thm:nat}, gives a natural transformation fitting in the square
\[
\xymatrix{
M_1G(n) \ar[r]^{e(n)} \ar[d]_{\sigma(n)} & G(n+1) \ar[d] \\
G \ar[r]_\eta & KG
}
\]
The existence of these natural isomorphisms will then be sufficient to conclude that the unit
$\eta:G \to KG$\/ is a weak equivalence.

Before going on to the construction we first recall some basic facts
about the kinds of models of type theory considered in this paper 
(the Hofmann-Streicher style groupoid models
and the combinatorial realizability models described in Section
\ref{subsec:CRM}).  These models are genuine denotational models in
the sense that each term, type and judgement is assigned a canonical
interpretation (as opposed to many realizability models where a given
judgement may have many different realizers).  A consequence of this fact, and of the
interpretation of substitution in these models, is that the
interpretation is \emph{compositional} in the sense that if we are
given an open judgement $\judge{x:A}{\mathcal{J}}$ and a term
$a:A$, then the interpretation of $\mathcal{J}[a/x]$ is completely
determined by the interpretations of $\judge{x:A}{\mathcal{J}}$ and
$a:A$.  In particular, in order to prove that
$\mathcal{J}[a/x]$ and $\mathcal{J}[b/x]$ receive the same
interpretation it suffices to show that $a$ and $b$ receive the same
interpretation.  E.g., in the Hofmann-Streicher style groupoid models
$\lscott b(a)\rscott$ is canonically determined as the
canonical section induced by $\lscott b(x)\rscott$ and
$\lscott a\rscott$ (cf.~the discussion in Section \ref{subsec:CRM}
following the proof of Theorem \ref{theorem:CRM}).  Explicitly, if
$\judge{\Gamma,x:A}{b(x):B(x)}$ and $\judge{\Gamma}{a:A}$, then we
have a pullback diagram
\begin{align*}
  \xy
  {\ar(0,15)*+{\lscott \Gamma, z:B(a)\rscott};(40,15)*+{\lscott
      \Gamma,x:A,z:B(x)\rscott}};
  {\ar(0,15)*+{\lscott
      \Gamma,z:B(a)\rscott};(0,0)*+{\lscott\Gamma\rscott}};
  {\ar_{\lscott\judge{\Gamma}{a:A}\rscott}(0,0)*+{\lscott\Gamma\rscott};(40,0)*+{\lscott\Gamma,x:A\rscott}};
  {\ar(40,15)*+{\lscott\Gamma,x:A,z:B(x)\rscott};(40,0)*+{\lscott\Gamma,x:A\rscott}};
  \endxy
\end{align*}
and a section $\lscott\judge{\Gamma,x:A}{b(x):B(x)}\rscott$ of the
projection
$\lscott\Gamma,x:A,z:B(x)\rscott\to\lscott\Gamma,x:A\rscott$.
This induces a canonical section $\lscott a\rscott^{*}(\lscott
b(x)\rscott)$ of the projection
$\lscott\Gamma,z:B(a)\rscott\to\lscott\Gamma\rscott$ which we
define to be the interpretation of $\judge{\Gamma}{b(a):B(a)}$.

Similar remarks apply to the combinatorial realizability models. We
will make use of this compositionality of the interpretations at
several places below.

We begin by constructing by induction a sequence of functors 
\[ \lscott-\rscott_n:M_1G(n) \to G \]
as well as a sequence of graph homomorphisms $\alpha(n):G(n) \to G$\/ such that
\[ 
\xymatrix{
M_1G(n) \ar[r]^-{e(n)} \ar[dr]_{\lscott-\rscott_n} & G(n+1) \ar[d]^{\alpha(n+1)}\\
& G
}
\]
is commutative.

For $n=-1$, we set $\alpha(-1):G(-1)=G \to G$\/ to be the identity. Note that
in order to specify the functor $\lscott-\rscott_n$\/ 
it suffices to specify a graph homomorphism $G(n) \to G$\/ and
then to use the Hofmann-Streicher semantics to extend this to $M_1G_n$. 
(Whence the notation overloading.) Thus at each stage, we may
let $\lscott-\rscott_n$\/ be the the functor induced by $\alpha(n)$. It remains to be shown then that
it coequalizes the relevant maps so that it factors through
the coequalizer $e(n)$\/ resulting in the desired $\alpha(n+1)$\/ as in the above diagram.

We begin with the base case, where we have to show that 
$\lscott-\rscott_{-1}:M_1G \to G$\/ makes the diagram
\begin{align*}
  \xy
  {\ar^{M_{1}\gamma}(0,0)*+{M_{1}\mathcal{F}G};(30,0)*+{M_{1}G}};
  {\ar_{M_{1}\Phi_{G}}@/_1em/(0,0)*+{M_{1}\mathcal{F}G};(15,-10)*+{M_{1}^{2}G}};
  {\ar_{\mu_{G}}@/_1em/(15,-10)*+{M_{1}^{2}G};(30,0)*+{M_{1}G}};
  {\ar^{\lscott-\rscott_{-1}}(30,0)*+{M_{1}G};(50,0)*+{G}};
  \endxy
\end{align*}
commute. Consider an object of $M_1\mathcal{F}G$, regarded as a term 
$\varphi(\namefor{g\star f})$, where $\star$\/ denotes the formal composition in the free groupoid.
 Here, we assume that $\varphi(x)$\/ is skeletal and that
$g,f$\/ are basic edges in $G$. (Technically we
  must consider arbitrary
  strings of formal composites in the free groupoid and a skeletal term $\varphi(x_1,\ldots, x_n)$, 
but the reasoning in that case is identical to the reasoning given here.)
On the one hand, $M_1\gamma$\/ sends this term to $\varphi(\namefor{g \circ f})$\/ 
(where $g \circ f$\/ is the composite in $G$\/ using the groupoid structure $\gamma$), while on the
other hand $\mu_G.M_1\Phi_G$\/ sends it to $\varphi(\namefor{g}\cdot\namefor{f})$\/ 
(where $\cdot$\/ is the formal composition in $M_1(G)$). 
Thus to prove that the two composites are equal, we must show that
\begin{align*}
  \lscott\varphi(\namefor{g\circ f})\rscott_{-1} 
   & = \lscott\varphi(\namefor{g}\cdot\namefor{f})\rscott_{-1}
\end{align*}
By compositionality of the interpretation it suffices to show that
$\lscott\namefor{g\circ f}\rscott_{-1}=\lscott\namefor{g}\cdot\namefor{f}\rscott_{-1}$,
which holds by functoriality of the interpretation.  

Next, assume that we have defined $\lscott-\rscott_n$\/ and $\alpha(n+1)$. We must then show that
$\lscott-\rscott_{n+1}$, obtained by interpreting $M_1G(n+1)$\/ in $G$, makes the diagram

\begin{equation}\label{eq:int_n}
  \xy
  {\ar@/^1em/^{\mu_{G(n)}}(0,0)*+{M_{1}^{2}G(n)};(15,10)*+{M_{1}G(n)}};
  {\ar@/^1em/^{M_{1}\eta_{G(n)}}(15,10)*+{M_{1}G(n)};(30,0)*+{M_{1}^{2}G(n)}};
  {\ar_{1_{M_{1}^{2}G(n)}}(0,0)*+{M_{1}^{2}G(n)};(30,0)*+{M_{1}^{2}G(n)}};
  {\ar^-{M_{1}e(n)}(30,0)*+{M_{1}^{2}G(n)};(60,0)*+{M_{1}G(n+1)}};
  {\ar^-{\lscott-\rscott_{n+1}}(60,0)*+{M_{1}G(n+1)};(90,0)*+{G}};
  \endxy
\end{equation}
commute.
To this end, we first establish the following useful lemma.

\begin{lemma}
The morphisms $\alpha(n+1)$\/ and the interpretations $\lscott-\rscott_{n+1}, \lscott-\rscott_{-1}$\/ 
form a commutative diagram
\begin{equation}\label{eq:alpha}
\xymatrix{
M_1G(n+1) \ar[rr]^{M_1\alpha(n+1)} \ar[drr]_{\lscott-\rscott_{n+1}} && M_1G \ar[d]^{\lscott - \rscott_{-1}} \\
&& G.
}
\end{equation}
\end{lemma}

\begin{proof}
By compositionality, it suffices to verify that the diagram commutes when we precompose with 
the unit $\eta:G(n+1) \to M_1G(n+1)$. But then by naturality of the unit and the definition of
$\lscott-\rscott_{n+1}$\/ we get
\[ \lscott-\rscott_{n+1}\circ \eta = \alpha(n+1) =\lscott-\rscott_{-1}\circ \eta \circ \alpha(n+1)=
\lscott-\rscott_{-1}\circ M_1(\alpha(n+1))\circ \eta
\]
as required. 
\end{proof}

Now consider the diagram
\[
\xymatrix{
M_1G(n) \ar[dr]_{\lscott -\rscott_n} \ar[dd]_{M_1\eta G(n)}  \ar[rr]^{M_1\alpha(n)} 
  && M_1G \ar[dl]_{\lscott -\rscott_{-1}} \\
& G \\
M_1^2G(n) \ar[rr]_{M_1e(n)} && M_1G(n+1) \ar[uu]_{M_1\alpha(n+1)} \ar[ul]^{\lscott -\rscott_{n+1}}
}
\]
The outer diagram commutes because by IH we have 
\[ \alpha(n+1) \circ e(n)\circ \eta=\alpha(n).\]
 The two triangles
commute by the argument given above. Therefore, to show that~\eqref{eq:int_n} commutes, it suffices
(again by compositionality)
to show that for a term $\namefor{s}$, where $s$\/ an element of $G(n)$, we have
$\lscott \mu\namefor{s}\rscott_{n}=\lscott M_1e(n)\namefor{s}\rscott_{n+1}$. But that is immediate
from $\mu\namefor{s}=s$\/ and the IH. 
This completes the proof that we have a well-defined sequence of functors
$\lscott-\rscott_n:M_1G(n) \to G$.

\begin{theorem}\label{thm:equiv}
  The adjunction $K\adjoint\Pi_{1}$ is a Quillen equivalence.
  \begin{proof}
    It suffices to show that $\eta_{G}$ is essentially surjective on
    objects and full. Denote the composite $M_1G(n) \to G(n+1) \to G$\/ by $\epsilon(n)$, and note
that this is actually a functor.
We may then consider, for each $n$, the (non-commutative) square
\[
\xymatrix{
M_1G(n) \ar[r]^{e(n)} \ar[d]_{\lscott-\rscott_n} & G(n+1) \ar[d] \\
G \ar[r]_\eta & KG.
}
\]
According to Theorem~\ref{thm:nat}, we may specify a natural transformation
in this square by giving the components at the basic terms of $M_1G(n)$\/ and verifying
that these are natural. We will do this first for $n=-1$. Then a basic term of $M_1G$\/ is simply an
element of $G$, and we take the component of the natural transformation to be the identity at that element.
This gives the natural transformation $\tau(-1):\eta\lscott-\rscott_{-1} \Rightarrow \epsilon(-1)$.

In the inductive step, we assume we have constructed $\tau(n):\eta\lscott-\rscott_{n} \Rightarrow
\epsilon(n)$, and we wish to define $\tau(n+1)$. Again by Theorem~\ref{thm:nat}, we only have to
specify the components at basic terms. Given such basic term $\namefor{[t]}$, where $t$\/ is an element
of $M_1G(n)$, we take this component to be $\tau(n)_t$. Naturality is then inherited from $\tau(n)$,
and it also is immediate that this is independent of the choice of representative of $[t]$\/ because
diagram~\eqref{eq:alpha} commutes.

    To see that $\eta_{G}$ is essentially surjective on objects, let
    $[t]$ in $KG$ be given.  So, $t$ is in some $G(n)$ and by the
    construction of the natural transformations above we get a component
    $\tau(n):[t] \to \eta\lscott-\rscott_n(t)$, where the latter is in the image of $\eta_{G}$, as required.
    Similarly, given an arrow $[g]:\eta(a)\to\eta(b)$ in $KG$\/
    we it follows that $g$ is in some $G(n)$. Because the components of $\eta(a)$\/ and $\eta(b)$\/ 
   of the natural transformation $\tau(n)$\/ are identities, 
the naturality square at $g$\/ of $\tau(n)$\/ then simply exhibits $g$\/
    as equal to a map in the image of
    $\eta_{G}$.
  \end{proof}
\end{theorem}


\appendix

\section{Rules of type theory}\label{sec:types_ap}

In this appendix we describe the syntax of the system $\mathbb{T}_{\omega}$.  All rules below are
stated in an ambient context which is omitted for ease of presentation.

\subsection{Structural rules}

\begin{prooftree}
  \AxiomC{$\judge{\Gamma}{\mathcal{J}}$}
  \RightLabel{Weakening}
  \UnaryInfC{$\judge{\Delta,\Gamma}{\mathcal{J}}$}
\end{prooftree}
where $\mathcal{J}$ ranges over judgements and we assume without loss
of generality that the variables declared in $\Delta$ and $\Gamma$ are
disjoint.

\begin{prooftree}
  \AxiomC{$a:A$}
  \AxiomC{$\judge{x:A,\Delta}{B(x):\type}$}
  \RightLabel{Type substitution}
  \BinaryInfC{$\judge{\Delta[a/x]}{B(a):\type}$}
\end{prooftree}
\begin{prooftree}
  \AxiomC{$a:A$}
  \AxiomC{$\judge{x:A,\Delta}{b(x):B(x)}{}$}
  \RightLabel{Term substitution}
  \BinaryInfC{$\judge{\Delta[a/x]}{b(a):B(a)}$}
\end{prooftree}

\begin{prooftree}
  \AxiomC{$A:\type$}
  \RightLabel{Variable declaration}
  \UnaryInfC{$\judge{x:A,\Delta}{x:A}$}
\end{prooftree}

\subsection{Rules governing definitional equality}

\vspace{.25cm}
\begin{center}
  \AxiomC{$A:\type$}
  \UnaryInfC{$A=A:\type$}
  \DisplayProof
\hspace{1cm}
\AxiomC{$A=B:\type$}
  \UnaryInfC{$B=A:\type$}
  \DisplayProof
\end{center}
\vspace{.5cm}
\begin{prooftree}
  \AxiomC{$A=B:\type$}
  \AxiomC{$B=C:\type$}
  \BinaryInfC{$A=C:\type$}
\end{prooftree}
\vspace{.5cm}
\begin{center}
  \AxiomC{$a:A$}
  \UnaryInfC{$a=a:A$}
  \DisplayProof
  \hspace{1cm}
  \AxiomC{$a=b:A$}
  \UnaryInfC{$b=a:A$}
  \DisplayProof
\end{center}
\vspace{.5cm}
\begin{prooftree}
  \AxiomC{$a=b:A$}
  \AxiomC{$b=c:A$}
  \BinaryInfC{$a=c:A$}
\end{prooftree}
\vspace{.5cm}
\begin{prooftree}
  \AxiomC{$a=b:A$}
  \AxiomC{$\judge{x:A}{B(x):\type}$}
  \BinaryInfC{$B(a)=B(b):\type$}
\end{prooftree}
\vspace{.5cm}
\begin{prooftree}
  \AxiomC{$a=b:A$}
  \AxiomC{$\judge{x:A}{f(x):B(x)}$}
  \BinaryInfC{$f(a)=f(b):B(a)$}
\end{prooftree}
\begin{prooftree}
  \AxiomC{$A=B:\type$}
  \AxiomC{$a:A$}
  \BinaryInfC{$a:B$}
\end{prooftree}

\subsection{Formation rules}

\vspace{.25cm}
\begin{prooftree}
  \AxiomC{$\judge{x:A}{B(x):\type}$}
  \RightLabel{$\prod$ formation}
  \UnaryInfC{$\prod_{x:A}B(x):\type$}
\end{prooftree}
\vspace{.5cm}
\begin{prooftree}
  \AxiomC{$\judge{x:A}{B(x):\type}$}
  \RightLabel{$\sum$ formation}
  \UnaryInfC{$\sum_{x:A}B(x):\type$}
\end{prooftree}
\vspace{.5cm}
\begin{prooftree}
  \AxiomC{$a,b:A$}
  \RightLabel{$\id{}$ formation}
  \UnaryInfC{$\judge{}{\idn{A}{}(a,b):\type}$}
\end{prooftree}
\vspace{.5cm}
\begin{prooftree}
  \AxiomC{}
  \RightLabel{$\NN$ formation}
  \UnaryInfC{$\judge{}{\NN:\type}$}
\end{prooftree}

\subsection{Introduction and elimination rules for dependent products}

\vspace{.25cm}
\begin{prooftree}
  \AxiomC{$\judge{x:A}{f(x):B(x)}$}
  \RightLabel{$\prod$ introduction}
  \UnaryInfC{$\lambda_{x:A}f(x):\prod_{x:A}B(x)$}
\end{prooftree}
\vspace{.5cm}
\begin{prooftree}
  \AxiomC{$f:\prod_{x:A}B(x)$}
  \AxiomC{$a:A$}
  \RightLabel{$\prod$ elimination}
  \BinaryInfC{$\app(f,a):B(a).$}
\end{prooftree}

\subsection{Introduction and elimination rules for dependent sums}

\vspace{.25cm}
\begin{prooftree}
  \AxiomC{$a:A$}
  \AxiomC{$b:B(a)$}
  \RightLabel{$\sum$ introduction}
  \BinaryInfC{$\pair(a,b):\sum_{x:A}B(x)$}
\end{prooftree}
\vspace{.5cm}
\begin{prooftree}
  \AxiomC{$\judge{}{p:\sum_{x:A}B(x)}$}
  \AxiomC{$\judge{x:A,y:B(x)}{\psi(x,y):C\bigl(\pair(x,y)\bigr)}$}
  \RightLabel{$\sum$ elimination}
  \BinaryInfC{$\RSigma\bigl([x:A,y:B(x)]\psi(x,y),p\bigr):C(p)$}
\end{prooftree}

\subsection{Introduction and elimination rules for identity types}

\vspace{.25cm}
\begin{prooftree}
  \AxiomC{$a:A$}
  \RightLabel{$\id{}$ introduction}
  \UnaryInfC{$r(a):\idn{A}{}(a,a)$}
\end{prooftree}
\vspace{.5cm}
\begin{prooftree}
  \AxiomC{$\judge{x:A,y:A,z:\idn{A}{}(x,y)}{B(x,y,z):\type}$}
  \noLine
  \UnaryInfC{$\judge{x:A}{\varphi(x):B\bigl(x,x,r(x)\bigr)}$}
  \noLine
  \UnaryInfC{$f:\idn{A}{}(a,b)$}
  \RightLabel{$\id{}$ elimination}
  \UnaryInfC{$\Jelim{[x,y:A,z:\idn{A}{}(x,y)]B(x,y,z)}\bigl([x:A]\varphi(x),a,b,f):B(a,b,f)$}
\end{prooftree}

\subsection{Introduction and elimination rules for natural numbers}

\vspace{.25cm}
\begin{prooftree}
  \AxiomC{}
  \RightLabel{$\NN$ introduction (i)}
  \UnaryInfC{$\zero:\NN$}
\end{prooftree}
\vspace{.5cm}
\begin{prooftree}
  \AxiomC{$n:\NN$}
  \RightLabel{$\NN$ introduction (ii)}
  \UnaryInfC{$\successor(n):\NN$}
\end{prooftree}
\vspace{.5cm}
\begin{prooftree}
  \AxiomC{$n:\NN$}
  \AxiomC{$c:C(\zero)$}
  \AxiomC{$\judge{x:\NN,y:C(x)}{\gamma(x,y):C\bigl(\successor(x)\bigr)}$}
  \RightLabel{$\NN$ elimination}
  \TrinaryInfC{$\rec\bigl(n,c,[x:\NN,y:C(x)]\gamma(x,y)\bigr):C(n)$}
\end{prooftree}

\subsection{Conversion rules}

\vspace{.25cm}
\begin{prooftree}
  \AxiomC{$\lambda_{x:A}f(x):\prod_{x:A}B(x)$}
  \AxiomC{$a:A$}
  \RightLabel{$\prod$ conversion}
  \BinaryInfC{$\app\bigl(\lambda_{x:A}f(x),a\bigr) \;=\; f(a):B(a)$}
\end{prooftree}
\vspace{.5cm}
\begin{prooftree}
  \AxiomC{$a:A$}
  \AxiomC{$b:B(a)$}
  \AxiomC{$\judge{x:A,y:B(x)}{\psi(x,y):C\bigl(\pair(x,y)\bigr)}$}
  \RightLabel{$\sum$ conversion}
  \TrinaryInfC{$\RSigma\bigl([x:A,y:B(x)]\psi(x,y),\pair(a,b)\bigr) \;=\; \psi(a,b):C\bigl(\pair(a,b)\bigr)$}
\end{prooftree}
\vspace{.5cm}
\begin{prooftree}
  \AxiomC{$a:A$}
  \RightLabel{$\id{}$ conversion}
  \UnaryInfC{$\Jelim{[x,y:A,z:\idn{A}{}(x,y)]B(x,y,z)}\bigl([x:A]\varphi(x),a,a,r(a)\bigr)\;=\;\varphi(a):B\bigl(a,a,r(a)\bigr)$}
\end{prooftree}
\vspace{.5cm}
\begin{prooftree}
  \AxiomC{}
  \RightLabel{$\NN$ conversion (i)}
  \UnaryInfC{$\rec\bigl(\zero,c,[x:\NN,y:C(x)]\gamma(x,y)\bigr) = c :C(\zero)$}
\end{prooftree}
\vspace{.5cm}
\begin{prooftree}
  \AxiomC{$n:\NN$}
  \RightLabel{$\NN$ conversion (ii)}
  \UnaryInfC{$\rec\bigl(\successor(n),c,[x:\NN,y:C(x)]\gamma(x,y)\bigr) = \gamma\bigl(n,\rec\bigl(c,[x:\NN,y:C(x)]\gamma(x,y),n\bigr)\bigr) :C\bigl(\successor(n)\bigr)$}
\end{prooftree}

\newcommand{\SortNoop}[1]{}
\providecommand{\bysame}{\leavevmode\hbox to3em{\hrulefill}\thinspace}
\providecommand{\MR}{\relax\ifhmode\unskip\space\fi MR }
\providecommand{\MRhref}[2]{%
  \href{http://www.ams.org/mathscinet-getitem?mr=#1}{#2}
}
\providecommand{\href}[2]{#2}

\end{document}